\documentclass[11pt,a4paper]{amsart}
\usepackage{amsmath,amsfonts,amssymb,epsf}
\usepackage[T1]{fontenc}

\usepackage{graphicx}

\let\bbbibitem\bibitem
\renewcommand{\bibitem}[2][]{\bbbibitem[#1]{#2}\label{#2}}

\def\fin{\hfill\hbox{\hskip .2cm $\square$}\medskip}

\newtheorem{theo}{Theorem}[section]
\newtheorem{lemma}[theo]{Lemma}
\newtheorem{prop}[theo]{Proposition}
\newtheorem{coro}[theo]{Corollary}
\newtheorem{defi}[theo]{Definition}
\newtheorem{propdef}[theo]{Proposition-definition}
\newtheorem{claim}[theo]{Claim}
\newenvironment{demo}[1][\hspace{-3pt}]{{\noindent\em Proof #1.~ }}{\fin}
\newenvironment{rema}{\smallskip\noindent{\em Remark.}}{\medskip}
\newenvironment{thm}[1][\hspace{-3pt}]{{\vspace{0.4cm}\noindent\bf 
Theorem #1}~~\em}{\vspace{0.4cm}}
\newenvironment{cor}[1][\hspace{-3pt}]{{\vspace{0.4cm}\noindent\bf 
Corollary #1}~~\em}{\vspace{0.4cm}}

\def\cal{\mathcal}
\def\a{\alpha}
\def\b{\beta}
\def\g{\gamma}
\def\R{{\mathbb R}} 
\def\d{{\rm d}}
\def\PU{{\rm PU}}
\def\SU{{\rm SU}}
\def\C{{\mathbb C}}
\def\B{{\mathbb H}_{\mathbb C}}
\def\Bm{{\mathbb H}_{\mathbb C}^m}    
\def\Bn{{\mathbb H}_{\mathbb C}^n}
\def\Bl{{\mathbb H}_{\mathbb C}^l}
\def\N{{\mathbb N}}
\def\Z{{\mathbb Z}}
\def\R{{\mathbb R}}
\def\fd{\longrightarrow}
\def\pfd{\rightarrow}

\def\la{\langle}
\def\ra{\rangle}
\def\ov{\overline}

\def\duzf{\partial^{1,0}f}
\def\dzuf{\partial^{0,1}f}
\def\duzfb{\partial^{1,0}f(z_\b)}

\def\dzufb{\partial^{0,1}f(z_\b)}
\def\om{\omega}
\def\za{{z_\alpha}}

\def\nab{\nabla}
\def\zb{{z_\beta}}
\def\zab{{\bar z_\alpha}}
\def\zbb{{\bar z_\beta}}
\def\df{{\rm d}f}
\def\Qrond{\stackrel{\,\circ}{Q}}
\def\gn{g_n}
\def\gm{g_m}
\def\t{\theta}
\def\ds{\displaystyle}
\def\e{\varepsilon}

\def\wt{\widetilde}
\def\G{\Gamma}

\renewcommand\Re{{\rm Re}}
\renewcommand\Im{{\rm Im}}

\title{Harmonic maps and representations of non-uniform lattices of $\PU(m,1)$}

\begin{document}

\author{Vincent Koziarz} 
\author{Julien Maubon}
\address{Institut Elie Cartan, Université
    Henri Poincaré, B. P. 239, F-54506 Vand\oe uvre-lès-Nancy Cedex,
    France}
\email{koziarz@iecn.u-nancy.fr}  
\email{maubon@iecn.u-nancy.fr}

\date{March 2004}


\sloppy


\begin{abstract}
  We study representations of lattices of $\PU(m,1)$ into $\PU(n,1)$. We
  show that if a representation $\rho$ is reductive and if $m\geq 2$, then
  there exists a finite energy harmonic $\rho$-equivariant map from 
  $\Bm$ to $\Bn$. This allows us to give a differential geometric proof of
  rigidity results obtained by M.~Burger and A.~Iozzi. We also define a new
  invariant associated to representations into $\PU(n,1)$ of non-uniform
  lattices of $\PU(1,1)$, and more generally of fundamental groups of
  orientable surfaces of finite topological type and negative Euler
  characteristic. We prove that this invariant is bounded by a constant
  depending only on the Euler characteristic of the surface and we give a
  complete characterization of representations with maximal invariant, thus
  generalizing the results of D.~Toledo for uniform lattices.       
\end{abstract}

\maketitle


\setcounter{section}{-1}

\section{Introduction}

Lattices in semi-simple Lie groups with no compact factor (say, defined
over $\R$ and with trivial center) enjoy several rigidity properties. For
example, with the exception of lattices in groups locally isomorphic to
${\rm PSL}(2,\R)$, they all satisfy Mostow strong rigidity, which roughly
means the 
following. Take two such Lie groups $G$ and $H$, an irreducible lattice
$\G$ in $G$, and a representation (that is, a homomorphism of groups) of
$\G$ into $H$. Assume that the representation is faithful and that the image
of $\G$ is also a lattice in $H$. Then the representation extends to a
homomorphism from the ambient Lie group $G$ to $H$ (see~\cite{Mo73}).  
Another type of rigidity, known as Margulis superrigidity, provides the
same kind of conclusion but with much weaker assumptions: the only
hypothesis is that the image of $\G$ should be Zariski-dense in
$H$. Superrigidity holds for lattices in Lie groups of rank at least 2
(\cite{Ma91}) and for lattices of quaternionic or octonionic hyperbolic
spaces (that is, lattices in the rank one Lie groups ${\rm Sp}(m,1)$,
$m\geq 2$, and ${\rm F}_4^{-20}$) (see~\cite{Corlette92} and~\cite{GS92}).
On the contrary, for lattices of real and complex hyperbolic spaces,
namely, lattices in the other rank one Lie groups ${\rm PO}(m,1)$ and 
${\rm PU}(m,1)$, superrigidity is known to fail in general.    

\medskip

In this paper, we will focus on lattices in ${\rm PU}(m,1)$, the group of
orientation-preserving isometries (or equivalently, of biholomorphisms) of
complex hyperbolic $m$-space $\Bm=\PU(m,1)/{\rm U}(m)$. They are of
particular interest  because they lie somewhere in between the very
flexible lattices of ${\rm PO}(m,1)$ and those, superrigid, of the higher
rank Lie groups.  

In~\cite{GM87}, W.~M.~Goldman and J.~J.~Millson studied representation
spaces of {\it uniform} torsion-free lattices $\G<{\rm SU}(m,1)$ (which can
be considered as lattices in $\PU(m,1)$) into ${\rm PU}(n,1)$, for $n>m\geq
2$. They proved that there are no non-trivial deformations of the standard
representation of such a lattice. This means that all nearby
representations are {\it $\C$-Fuchsian}, namely, they are discrete,
faithful, and they stabilize a totally geodesic copy of $\Bm$ in $\Bn$. The
case $m=1$ was previously treated by Goldman in~\cite{Go85}.
Note that the corresponding statement for lattices in ${\rm PO}(m,1)$ is false
(cf. for example~\cite{JM87}).  

They also conjectured that a much stronger rigidity should hold. The
{\it volume} of a representation $\rho$ of a torsion-free uniform lattice  
$\G<{\rm PU}(m,1)$ into  ${\rm PU}(n,1)$ is defined by pulling-back the
Kähler form of $\Bn$ on $\Bm$ via the representation, taking its $m$-th
exterior power to obtain a de Rham cohomology class in
$H^{2m}_{DR}(\G\backslash\Bm)$  and evaluating it on the fundamental class
of the compact quotient $\G\backslash\Bm$. Observe that if $\G<{\rm SU}(m,1)$
and if $\rho:\G\fd\PU(n,1)$, $n>m$, is the standard representation, then 
${\rm vol}(\rho)={\rm Vol}(\G\backslash\Bm)$. Their conjecture then
reads: any representation $\rho$ such that ${\rm vol}(\rho)={\rm
  Vol}(\G\backslash\Bm)$ must be $\C$-Fuchsian. This was proved
by K.~Corlette in~\cite{Co88} for $m\geq 2$ and by D.~Toledo in~\cite{To89}
for $m=1$. Remark that the volume assumption is needed precisely because
lattices in $\PU(m,1)$ are not superrigid. 

Recently, M. Burger and A. Iozzi proved in~\cite{BI01} (see
also~\cite{Io02}) that the conjecture 
is also true for {\it non-uniform} lattices of $\PU(m,1)$, $m\geq 2$, if
one suitably modifies the definition of the ``volume'' of the
representation (indeed, with the former one, any representation of a
non-uniform lattice has zero volume). We will explain precisely how this
invariant is computed in section~\ref{invariant} but here we sketch its
definition. Again, the Kähler form $\om_n$ of $\Bn$ is pulled-back to the
quotient $\G\backslash\Bm$ via the representation. It turns out that this
gives a well-defined $L^2$-cohomology class in $H^{2}_{(2)}(\G\backslash\Bm)$.
Now, integrating a $L^2$-representative
$\rho^\star\om_n$ against the Kähler form $\om_m$ of $\G\backslash\Bm$,
we get the Burger-Iozzi invariant (slightly modified from~\cite{BI01}):
$$
\tau(\rho):=\frac{1}{2m}\int_{\G\backslash\Bm} \la\rho^\star\om_n,\om_m\ra dV_m~.
$$      
 In complex dimension~1 and for uniform lattices, this invariant coincides
 with the invariant defined in~\cite{To89}. 
We can now state the main theorem of~\cite{BI01}:

\begin{thm}[A] Let $\G$ be a torsion-free lattice in $\PU(m,1)$, $m\geq 2$,
  and let $\rho:\G\fd\PU(n,1)$ be a representation. Then $|\tau(\rho)|\leq
  {\rm Vol}(\G\backslash\Bm)$ and equality holds if and only if there
  exists a totally geodesic isometric $\rho$-equivariant embedding of $\Bm$
  into $\Bn$ (in particular, $\rho(\G)$ seen as a subgroup of $\PU(m,1)$ is
  a lattice).
\end{thm}

Burger and Iozzi's proof heavily relies on the theory of
bounded cohomology developed by Burger and N.~Monod in~\cite{BM02}.   
As a corollary, they obtain the result of Goldman and
Millson for a general lattice:

\begin{cor}[A']
  Let $\G$ be a torsion-free lattice in $\SU(m,1)$, $m\geq 2$, and let
  $n>m$. Then there are no non-trivial deformations of the standard
  representation of $\G$ into $\PU(n,1)$.
\end{cor}

The aim of this paper is to use harmonic map techniques to give a new and
more (differential) geometric proof of Theorem A and to extend this result
to the case of complex dimension 1, that is, of non-uniform lattices of
$\PU(1,1)$.    

The over-all harmonic map strategy for proving rigidity results about 
representations of lattices in a Lie Group $G$ to another Lie group
$H$ goes as follows. First, one has to know that there exists a
harmonic map between the corresponding  symmetric spaces,
equivariant w.r.t. the representation. Then, one must prove, generally by
using a Bochner-type formula, that there are additional constraints on the
harmonic map, which force it to be pluriharmonic, holomorphic, totally
geodesic, or isometric, depending on the situation. 

For a uniform $\G$ and when the target
symmetric space is non-positively curved (which will be assumed from now
on), the existence results for harmonic maps go back to J.~Eells
and J.~H.~Sampson in~\cite{ES64} and have been extended by several authors, in
particular by Corlette in~\cite{Co88}. The second step was pionneered by
Y.-T.~Siu in~\cite{Si80} where he proved a strenghtened version of Mostow
strong rigidity theorem in the case of Hermitian locally symmetric
spaces. This has later on been applied in different directions
by many authors. We should mention the proof of the above conjecture of
Goldman and Millson by Corlette in~\cite{Co88} and the geometric proof of
Margulis superrigidity theorem in the Archimedean setting worked out by
N.~Mok, Y.-T.~Siu and S.-K.~Yeung in~\cite{MSY93}.   

When the lattice is not uniform, the only general existence theorem for
harmonic maps is due to Corlette in~\cite{Corlette92}, see
Theorem~\ref{existharmonic} below. The main
issue is that to apply this theorem, one needs to prove
that there exists an equivariant map of {\it finite energy} (see
section~\ref{energiefinie} for the definition). If this is the case, the
harmonic map also has finite energy and the second step generally goes as
if the lattice was uniform, but is technically more involved. The energy
finiteness condition is very important, and in general difficult to
prove. 
In some particular cases it is possible to obtain a harmonic map by
other means (see for example~\cite{JZ97} and section~\ref{dim1} of this
paper) but then its energy is infinite and  
the analysis that follows becomes much harder. These are the reasons why, for 
example, ``geometric superrigidity'' for non-uniform lattices is not yet
proved. 

\vspace{.2cm}    
     
Our paper is organized as follows. The first three sections are devoted to
the proof of Theorem A. In section~\ref{energiefinie} we give
the necessary definitions and we prove that Corlette's general theorem
applies in our setting, so that we obtain our main existence theorem (cf. 
Theorem~\ref{finiteenergy}):

\begin{thm}[B] Let $\G$ be a torsion-free lattice in $\PU(m,1)$, $m\geq 2$,
  and $\rho:\G\fd\PU(n,1)$ be a representation such that $\rho(\Gamma)$ has  
  no fixed point on the boundary at infinity of $\Bn$. Then there exists a
  finite energy harmonic $\rho$-equivariant map from $\Bm$ to $\Bn$.   
\end{thm}

In section~\ref{pluri} we prove that the harmonic map previously obtained
is pluriharmonic and even holomorphic or antiholomorphic if its rank is high
enough (at least 3 at some point). Section~\ref{rigid} is devoted to the
precise definition of the Burger-Iozzi invariant and to the proof of
Theorem~A.

In section~\ref{dim1}, we study the case of lattices of
$\PU(1,1)$, that is, of fundamental groups of Riemann surfaces
with  a finite volume hyperbolic metric. The analogue of Theorem~A for {\it
  uniform} lattices was proved by Toledo in~\cite{To89}. In~\cite{GP03}
(see also~\cite{GP00}), N.~Gusevskii and J.~R.~Parker claim 
that if one restricts to {\it type-preserving} representations, then the
original definition of the Toledo invariant can be used to generalize
Toledo's result to non-uniform lattices. However, it seems to us that this 
claim is not entirely exact (see for example the remark following
Proposition~\ref{invharmonic}).
 
There are mainly two reasons why the 1-dimensional case is 
different from the higher dimensional one. First of all, Toledo and/or
Burger-Iozzi invariants are not defined for non-uniform lattices. Secondly,
there are representations for which no equivariant map of {\it finite
  energy} exists. It should also be noted that Corollary~A' fails in this
case by a result of Gusevskii and Parker (cf.~\cite{GP00}). 

As we shall see, it is in fact more natural to work in the general setting of
fundamental 
groups of orientable surfaces of finite topological type, that is surfaces
obtained by removing finitely many points from closed 
orientable surfaces. Using cohomology with compact support, we define at the
beginning of section~\ref{dim1} a new
invariant associated to representations of these fundamental groups into
$\PU(n,1)$, which we again call $\tau$. We obtain (see Theorem~\ref{surfaces}):

\begin{thm}[C] Let $\G$ be the fundamental group of a $p$-times punctured
  closed orientable surface $M$ of negative Euler characterisic $\chi(M)$, 
  and let $\rho:\G\fd\PU(n,1)$ be a representation. Then $|\tau(\rho)|\leq
  -2\pi\chi(M)$ and equality holds if and only if $\rho(\G)$ stabilizes
  a complex geodesic $L$ in $\Bn$, $\rho$ is faithful and discrete, and $M$
  is diffeomorphic to the quotient $\rho(\G)\backslash L$.
\end{thm}  

The proof relies on the fact that though there may be no equivariant
map of finite energy, there exists an equivariant harmonic map
whose energy density can be controlled. This control allows us to extend
the proofs given in the finite energy case to this setting. 


\begin{rema} In an earlier version of this paper, Theorem~C was proven in
a weaker form, and only for what we call {\it tame} representations (see
Definition~\ref{deftame}). M.~Burger and A.~Iozzi then
informed us that their methods should allow them to get rid of this
tameness assumption. Later, they communicated us the text~\cite{BI03}, where
they define a ``bounded Toledo number'' and prove Theorem C.  
\end{rema}

{\it \noindent Acknowledgments. } We would like to thank J.-P.~Otal who suggested
that it would be interesting to have a more geometric proof of Burger and
Iozzi's result. We also are grateful to F.~Campana and J.~Souto for helpful
conversations. We finally thank M.~Burger and A.~Iozzi for their interest
in our work and for having encouraged us to improve the first
draft of Theorem~C.

\section{Existence of finite energy equivariant harmonic maps}\label{energiefinie}

In this section, we assume that $m\ge 2$.

Let $\G$ be a torsion-free lattice in ${\rm PU}(m,1)$, the group of
biholomorphisms of complex hyperbolic $m$-space $\Bm$ and let
$\rho:\G\fd {\rm PU}(n,1)$ be a representation into the group of
biholomorphisms of complex hyperbolic $n$-space $\Bn$.  

We call $M$ the quotient manifold $\G\backslash\Bm$.
The representation $\rho$ determines a flat bundle
$M\times_\rho\Bn$ over $M$ with fibers isomorphic to $\Bn$. Since $\Bn$ is
contractible, this bundle has global sections. This is equivalent to the
existence of  maps (belonging to the same homotopy class) from $\Bm$ to
$\Bn$, equivariant w.r.t. the representation $\rho$. Let $f$ be such a map
(or section). 

We can consider the differential $\df$ of $f$ as a  
$f^\star T\Bn$-valued 1-form on $\Bm$. There is a natural pointwise scalar
product on 
such forms coming from the Riemannian metrics $g_m$ and $g_n$ (of constant
holomorphic sectional curvature $-1$) on $\Bm$ and
$\Bm$: if $(e_i)_{1\leq i\leq 2m}$ is a $g_m$-orthonormal basis of $T_x\Bm$,
then $\|\df\|^2_x:=\sum_i \gn(\df(e_i),\df(e_i))$. Since $f$ is
$\rho$-equivariant and the action of $\G$ on $\Bn$ via $\rho$ is isometric,
$\|\df\|$ is a well-defined function on $M$. We say that $f$ has {\it finite
energy} if the energy density $e(f):=\frac{1}{2}\|\df\|^2$ of $f$ is
integrable on $M$:
$$
E(f)=\frac{1}{2}\int_M \|\df\|^2 dV_m\,<\,+\infty~,
$$   
where $dV_m$ is the volume density of the metric $g_m$. When there is no risk of
confusion, we will write $e$ instead of $e(f)$ for the energy density of $f$. 

There is also a natural  connection
$\nab$ on $f^\star T\Bn$-valued 1-forms on $\Bm$  
coming from the Levi-Civita connections $\nab^m$ and
$\nab^n$ of $\Bm$ and $\Bn$. If $\nab^{f^\star T\Bn}$ denotes the
connection induced by $\nab^n$ on the bundle $f^\star T\Bn\fd\Bm$, then
$\nab\df(X,Y)=\nab^{f^\star T\Bn}_{X}\df(Y)-\df(\nab^m_X Y)$. Since $\nab^m$ and
$\nab^n$ are torsion-free, $\nab\df$ is a symmetric 2-tensor taking values
in $f^\star T\Bn$.    

A map $f:\Bm\fd\Bn$ is said to be {\it harmonic} if  ${\rm
  tr}_{g_m}\nab\df=0$.

The following theorem of Corlette~(\cite{Corlette92}) implies that if there
exists a finite energy $\rho$-equivariant map from the universal cover
$\Bm$ of $M$ to $\Bn$, and 
under a very mild assumption on $\rho$, then there exists a harmonic
$\rho$-equivariant map of finite energy from $\Bm$ to $\Bn$:

\begin{theo}\label{existharmonic}
  Let $X$ be a complete Riemannian manifold and $Y$ a complete
  simply-connected manifold with non-positive sectional curvature. Let
  $\rho:\pi_1(X)\fd {\rm Isom}(Y)$ be a representation such that
  the induced action of $\pi_1(X)$ on the sphere at infinity of $Y$ has no
  fixed point ($\rho$ is then called reductive). If there exists a
  $\rho$-equivariant map of finite energy 
  from the universal cover $\widetilde X$ of $X$ to $Y$, then there exists
  a harmonic $\rho$-equivariant map of finite energy from $\widetilde X$ to
  $Y$.  
\end{theo}

Theorem~B will therefore follow from the 

\begin{theo}\label{finiteenergy}
Let $\G$ be a torsion-free lattice in ${\rm PU}(m,1)$, $m\geq 2$, and let $\rho$
be a representation of $\G$ into ${\rm PU}(n,1)$. Then   
there exists a finite energy $\rho$-equivariant map $\Bm\fd\Bn$.
\end{theo}

\begin{demo}
Of course this is trivially true if the manifold is compact, that is, if
$\G$ is a uniform lattice. 
To prove the theorem in the non-uniform case, we recall some known
facts about the structure at infinity of the finite volume complex
hyperbolic manifold $M=\G\backslash\Bm$, cf. for example~\cite{Goldman},
or~\cite{Biquard} and~\cite{Hummel-Schroeder}.  

We will work with the Siegel model of complex hyperbolic space:
$$
\Bm=\{(z,w)\in\C^{m-1}\times\C\,\,\mid\,\,
2\Re(w)>\la\la z,z\ra\ra\}~,   
$$
where $\la\la.,.\ra\ra$ is the standard Hermitian product on $\C^{m-1}$. We
will call $h$ 
the function given by $h(z,w)=2\Re(w)-\la\la z,z\ra\ra$. The
boundary at infinity of $\Bm$ 
is the set $\{h=0\}\cup\{\infty\}$ and the horospheres in 
$\Bm$ centered at $\infty$  are the level sets of $h$.  
The complex hyperbolic metric (of constant holomorphic sectional curvature
$-1$) in the Siegel model of $\Bm$ is given by  
$$
g_m=\frac{4}{h(z,w)^2}\Big[(dw-\la\la dz,z\ra\ra)(d\bar{w}-\la\la
z,dz\ra\ra)+h(z,w)\la\la dz,dz\ra\ra\Big]~.
$$  

The stabilizer $P$ of 
$\infty$ in ${\rm PU}(m,1)$ is the 
semi-direct product ${\cal N}^{2m-1}\rtimes ({\rm U}(m-1)\times \{\phi_s\}_{s\in\R})$
where ${\cal N}^{2m-1}$ is 
the $(2m-1)$-dimensional Heisenberg group, ${\rm U}(m-1)$ is the unitary group and
$\{\phi_s\}_{s\in\R}$ is the one-parameter group corresponding to the
horocyclic flow associated to $\infty$. 
The group ${\cal N}^{2m-1}$ is a central extension of $\C^{m-1}$ and can be
seen as $\C^{m-1}\times\R$ with product given by
$(\xi_1,\nu_1)(\xi_2,\nu_2)=(\xi_1+\xi_2,
\nu_1+\nu_2+2\Im\la\la\xi_1,\xi_2\ra\ra)$. This is a two-step
nilpotent group which acts simply transitively and isometrically on horospheres.  Its
center ${\cal Z}$ is the group of ``vertical translations'': $\{(0,\nu),\,\nu\in\R\}$.

If we set $u+iv=2\ov{w}-\la\la z,z\ra\ra$, we obtain the
so-called {\it horospherical coordinates} $(z,v,u)\in
\C^{m-1}\times\R\times\R_+^\star$, in which the action of $P$ on $\Bm$ 
is given by:  
$$
(\xi,\nu)A\phi_s.(z,v,u)
=(Ae^{-s}z+\xi,e^{-2s}v+\nu+2\Im\la\la\xi,Ae^{-s}z\ra\ra,e^{-2s}u)
$$
and the metric $g_m$ takes the form
$$
g_m=\frac{du^2}{u^2}+\frac{1}{u^2}{\Big(-dv+2\Im\la\la z,dz\ra\ra\Big)}^2
+\frac{4}{u}\la\la dz,dz \ra\ra~.
$$
Replacing $u$ by $t=\log u$, the metric tensor decomposes as:
$$
g_m=dt^2+e^{-2t}{\Big(-dv+2\Im\la\la z,dz\ra\ra\Big)}^2
+4e^{-t}\la\la dz,dz \ra\ra~.
$$
The coordinates $(z,v,t)\in\C^{m-1}\times\R\times\R$ will also be called
horospherical coordinates.

A complex hyperbolic manifold $M$ of finite volume is the union of a
compact part and a finite number of disjoint cusps. 
Each cusp $C$ of $M$ is diffeomorphic to the product $N\times[0,+\infty)$, where
$N$ is a compact quotient of some horosphere $HS$ in $\Bm$. We can assume that $HS$ is
centered at ${\infty}$. The fundamental
group $\G_C$ of $C$, hence of $N$, can be identified with the stabilizer in
$\G$ of the horosphere $HS$: it is therefore equal to $\G\cap ({{\cal
  N}^{2m-1}\rtimes {\rm U}(m-1)})$. 

If we call $\beta$ the 1-form $-dv+2\Im\la\la z,dz\ra\ra$ on $\Bm$, it is
easily checked that ${\rm d}^c t:=J{\rm d}t=e^{-t}\beta$. Therefore, since ${\cal N}^{2m-1}\rtimes
{\rm U}(m-1)$ preserves the horospheres, $t$, $dt^2$, and $\beta$ are invariant by
$\G_C$. The decomposition of $g_m$ hence goes down to the cusp $C$ and we
have:  
$$
g_m = dt^2 + e^{-2t} \beta^2 + e^{-t} g ~,
$$ 
where $g$ is the image of $4\la\la dz,dz \ra\ra$.

\begin{rema} The Kähler form $\omega_m$, which we normalize so that
  $\omega_m(X,JX)\geq 0$, is of course exact on $\Bm$. More
precisely, $\omega_m=-{\rm d}{\rm d}^c t=-{\rm d}(e^{-t}\beta)$. The invariance
of $t$ and $\beta$ implies that this relation still holds in the cusps of $M$. 
\end{rema}

For lattices in ${\rm Sp}(m,1)$, $m\geq 2$, or in ${\rm F}_4^{-20}$,
Corlette proves in~\cite{Corlette92} a simple lemma that allows him to
deduce the existence of finite energy equivariant maps. Here, the same idea
will only provide the result for $m\geq 3$:

\begin{lemma}\label{retraction}
  Assume $m\geq 3$. Then there exists a finite energy retraction of
  $M=\G\backslash \Bm$ onto a compact subset of $M$. 
\end{lemma}

\begin{demo}
It is enough to construct the retraction on a cusp $C=N\times[0,+\infty)$ of
$M$: we define $r:N\times[0,+\infty)\fd N\times\{0\}$ obviously by
$r(x,t)=(x,0)$. 

If $(\frac{\partial}{\partial t},\frac{\partial}{\partial v},e_3,\ldots,e_{2m})$ is an
orthonormal basis of $T_{(x,0)}C$ compatible with the splitting of $\gm$, then 
 $(\frac{\partial}{\partial t},e^{t}\frac{\partial}{\partial
   v},e^{t/2}e_3,\ldots,e^{t/2}e_{2m})$ is such a basis of $T_{(x,t)}C$. Hence,
$
\|{\rm d}r\|^2_{(x,t)}=e^{2t}+(2m-2)e^{t}
$.
If we call $dV_N$ the volume element of $N\times\{0\}$, the volume element of
$N\times\{t\}$ is given by
$e^{\frac{1}{2}(-2t-(2m-2)t)}dV_N=e^{-mt}dV_N$. Hence, the energy of $r$ on
$C$ is  
$$
\frac{1}{2}\int_C \|{\rm d}r\|^2\,dV_m = \frac{1}{2}\int_0^{+\infty}\int_N
(2(m-1)e^{t}+e^{2t})e^{-mt}\, dV_N\,dt~.
$$
This is clearly finite if $m\geq 3$.      
\end{demo}

This retraction lifts to a map $\tilde r:\Bm\fd\Bm$, invariant by
$\G$. Therefore, if $f:\Bm\fd\Bn$ is any $\rho$-equivariant map, so
is $f\circ \tilde r$, and its energy is finite. The theorem is therefore proved if
$m\geq 3$.

\vspace{0.2cm}

In the case $m=2$, the energy density of the retraction $r$ 
grows like $e^{2t}$  when $t$ goes to infinity whereas the volume element
grows like $e^{-2t}$: the 
energy of $r$ is infinite and we need a deeper analysis of the situation in
the cusps.  

We fix a cusp $C=N\times[0,+\infty)$ of $M$ and we look for a finite
energy map from the universal cover $HS\times[0,+\infty)$ of $C$ to $\Bn$,
equivariant w.r.t. the fundamental group $\G_C$ of $C$ (equivalently, a
section of the restriction of the flat bundle $M\times_\rho\Bn$ to $C\subset M$). 

As we said, $\G_C$ can be seen as a subgroup of ${\cal N}\rtimes {\rm
  U}(1)$, where now ${\cal N}:={\cal N}^3$ is just $\C\times\R$. 
It follows
from L.~Auslander's generalization of Bieberbach's theorem (cf. \cite{Auslander})
that $\G_{\cal N}:=\G_C\cap{\cal N}$ is a discrete uniform subgroup
of ${\cal N}$, of finite index in $\G_C$. Therefore (\cite{Auslander}, Lemma
1.3.), $\G_{\cal N}$ cannot be contained in any proper analytic subgoup of
${\cal N}$. From this, it is easy to deduce that there exists $\varepsilon>0$
such that, for all $\g=(\xi_\g,\nu_\g)\in\G_{\cal N}$, $|\xi_\g|>\varepsilon$ as
soon as $\xi_\g\neq 0$. In other words, the image of the homomorphism of groups
$T:\,\G_{\cal N}\fd\C$, $\g\longmapsto\xi_\g$, is a lattice in $\C$. Let
$\g_1=(\xi_1,\nu_1)$ and $\g_2=(\xi_2,\nu_2)$ be two elements of $\G_{\cal N}$
such that $\xi_1$ and $\xi_2$ generate the lattice $T(\G_{\cal N})$.  
A straightforward computation yields
$[\g_1,\g_2]:=\g_1\g_2\g_1^{-1}\g_2^{-1}=(0,-2\Im(\xi_1\ov\xi_2))$. Since
$\xi_1$ and $\xi_2$ are linearly independent (over $\R$),
$\Im(\xi_1\ov\xi_2)\neq 0$ and hence the subgroup $\G_{\cal Z}:=\G_C\cap{\cal
  Z}$ of $\G_{\cal N}$ is non trivial. It is therefore isomorphic to ${\mathbb
  Z}$ and we call $\g_0$ its generator.  
Note that $\g_0$, $\g_1$ and $\g_2$ generate $\G_{\cal N}$.

The construction of the equivariant map will depend on the type of
$\rho(\g_0)$. Recall that an isometry of $\Bn$ can be of one of the
following (exclusive) 3 types: 
\begin{itemize}
\item {\it elliptic } if it has a fixed point in $\Bn$;
\item{\it parabolic} if it has exactly one fixed point on the sphere at
  infinity of $\Bn$ and no fixed points in $\Bn$;   
\item {\it hyperbolic} if it has exactly two fixed points on the sphere at
  infinity of $\Bn$ and no fixed points in $\Bn$. In this case, the isometry acts by
  translation on  the geodesic joining its fixed points at infinity. 
\end{itemize}

\begin{claim}\label{pashyperbolique}
  $\rho(\g_0)$ can not be a hyperbolic isometry of $\Bn$.
\end{claim}

\begin{demo}
  Assume that $\rho(\g_0)$ is a hyperbolic isometry of $\Bn$ and call $A_0$ 
  its axis (the geodesic joining its fixed points). Then, since $\g_1$ and
  $\g_2$ commute with $\g_0$, their images by $\rho$ commute with
  $\rho(\g_0)$, hence
  they must fix $A_0$ and act on it by translations: there exist
  $\tau_1,\tau_2\in\R$ such that $\rho(\g_1) A_0(t)=A_0(t+\tau_1)$ and
  $\rho(\g_2) A_0(t)=A_0(t+\tau_2)$. This implies that $\rho([\g_1,\g_2])$ acts
  trivialy on $A_0$. But $[\g_1,\g_2]={\g_0}^p$ for some $p\in{\mathbb Z}^\star$ and
  $\rho(\g_0)$ does not act trivially on $A_0$. This is a contradiction.  
\end{demo}

Hence $\rho(\g_0)$ is either elliptic or parabolic.
In both cases we will start
by constructing an equivariant map from the universal cover $HS\simeq {\cal
  N}$ of $N$ and then we shall extend it to the universal cover of the whole
cusp.   

\vspace{0.2cm}

{\it Case 1: $\rho(\g_0)$ is parabolic.} 
 The idea is to find an equivariant
map from $HS$ to a horosphere in $\Bn$ centered at the fixed point of
$\rho(\g_0)$ on the 
sphere at infinity $\partial\Bn$ of $\Bn$ and then
to extend it to $HS\times[0,+\infty)$ using the horocyclic flow defined by
the fixed point. Roughly speaking, when $t$ goes to infinity, the image of
$HS\times\{t\}$ must go to infinity in $\Bn$ fast enough so that the decay
of the metric in $\Bn$ prevents the energy density of the map from growing too
quickly. 

Using again the Siegel model for
$\Bn$, we may assume that the fixed point of $\rho(\g_0)$ is $\infty$. Since
$\g_0$ is in the center of $\G_C$, the whole group $\rho(\G_C)$ must fix
$\infty$, and therefore must be contained in its stabilizer in
${\rm PU}(n,1)$. Moreover, $\rho(\G_C)$ must stabilize each horosphere centered at
$\infty$. For, if this was not the case, there would be an element $\g\in\G_C$
such that $\rho(\g)$ is hyperbolic. But then, since $\g_0$ commutes with $\g$,
$\rho(\g_0)$ would fix the axis of $\rho(\g)$. This is impossible since we
assumed that $\rho(\g_0)$ is parabolic.  

We see $\Bn$ as $\C^{n-1}\times\R\times\R$ with horospherical coordinates
$(z',v',t'=\log u')$. The metric $g_n$ at a point
$(z',v',t')$ is given by $g_n=dt'^2+e^{-2t'}(-dv'+2\Im\la\la z',dz'\ra\ra)^2
+4e^{-t'}\la\la dz',dz' \ra\ra$.  

Let $HS'\subset\Bn$ be the horosphere $\C^{n-1}\times\R\times\{0\}$. The
representation $\rho$
can be seen as a homomorphism from the fundamental group of $N$ to the
isometry group of $HS'$. Since $HS'$ is contractible, there exists a
$\rho$-equivariant map $\varphi$ from the universal cover $HS\subset\Bm$ of
$N$ to the horosphere $HS'\subset\Bn$. 
Now, define a $\rho$-equivariant map $f$ from the universal cover
$HS\times[0,+\infty)$ of the cusp $C=N\times[0,+\infty)$ to $\Bn$ by: 
$$
\begin{array}{rcl}
f:HS\times[0,+\infty) & \fd & HS'\times[0,+\infty)\,\,\subset\Bn\\
(x,t) & \longmapsto & (\varphi(x),2t)
\end{array}
$$
Using the same notation as in Lemma \ref{retraction}, the energy density of $f$
can be estimated as follows:
$$
\begin{array}{rcl}
\|{\rm d}f\|^2_{(x,t)} & = & |{\rm d}f(\frac{\partial}{\partial
  t})|^2_{(\varphi(x),t)}
+e^{2t}|{\rm d}f(\frac{\partial}{\partial v})|^2_{(\varphi(x),t)}+
e^t\sum_{k=3}^{4}|{\rm d}f(e_k)|^2_{(\varphi(x),t)}\\
& \leq & 4+e^{-2t}\Big(e^{2t}|{\rm d}\varphi(\frac{\partial}{\partial
  v})|^2_{(\varphi(x),0)}+ e^t\sum_{k=3}^{4}|{\rm
  d}\varphi(e_k)|^2_{(\varphi(x),0)}\Big)\\ 
& \leq & 4+\|{\rm d}\varphi\|^2_x
\end{array}
$$  
where $\|{\rm d}\varphi\|$ denotes the norm of the differential of
$\varphi:HS\fd HS'$ computed with the metrics of $HS$ and $HS'$
induced from $g_m$ and $g_n$. 

The energy of $f$ in the cusp $C$ is therefore finite since:
$$
E_C(f)=\frac{1}{2}\int_C \|{\rm d}f\|^2 dV_m\leq
\frac{1}{2}\int_0^{+\infty}\int_N 
\Big(4+\|{\rm d}\varphi\|^2\Big)e^{-2t}dV_Ndt < +\infty~.
$$
 
\vspace{0.2cm}

{\it Case 2: $\rho(\g_0)$ is elliptic.} In this case, there is no canonical
``direction''  in which to send the slices $HS\times\{t\}$ to infinity
in $\Bn$.  Once the 
equivariant map $f$ is constructed on $HS\times\{0\}$, 
the most natural way to define it on $HS\times\{t\}$ is to set
$f_{|HS\times\{t\}}=f_{|HS\times\{0\}}$. Therefore, the growth of
the energy density in the cusp cannot be controled by some decay of the metric
in $\Bn$, and we must control it at the start. We shall achieve this by
demanding the equivariant map $HS\fd\Bn$ to be constant in the ``vertical
direction'' $\R$ of $HS=\C\times\R$. 

As mentionned before, $\G_{\cal N}$ is a finite index subgroup of $\G_C$ and we
have the tower of coverings:
$$
HS=\C\times\R\stackrel{\G_{\cal N}}{\longrightarrow}\widehat
N\stackrel{\G_C/\G_{\cal N}}{\longrightarrow}N~, 
$$
where $\widehat N=(\C\times\R)/{\G_{\cal N}}$ is a circle bundle over the
2-torus ${\mathbb T}=\C/T(\G_{\cal N})$, and $\G_C/\G_{\cal N}$ can be seen
as a finite subgroup of ${\rm U}(1)$, acting freely on this bundle. 

The group $\G_C/\G_{\cal N}$ is generated by a primitive $p$-th root of unity $a$ 
and its action on $\C$ preserves the lattice $T(\G_{\cal N})\subset\C$. This
implies that $a$ is a root of a degree 2 polynomial with integer coefficients
and hence the possible values of $a$ are $1$, $-1$, $e^{i\frac{2\pi}{3}}$, $i$,
or $e^{i\frac{\pi}{3}}$. On the other hand, the number of possible lattices is
also restricted: 
\begin{itemize}
\item if $a=1$ or $a=-1$, $T(\G_{\cal N})$ can be any lattice of $\C$;
\item if $a=i$, $T(\G_{\cal N})$ must be a square lattice, meaning that we can
  choose the first two generators $\g_1=(\xi_1,\nu_1)$ and $\g_2=(\xi_2,\nu_2)$ of
  $\G_{\cal N}$ so that $\xi_2=i\xi_1$; 
\item if $a=e^{i\frac{2\pi}{3}}$ or $a=e^{i\frac{\pi}{3}}$, $T(\G_{\cal
    N})$ must be an equilateral triangle lattice, meaning that we can
  choose $\g_1$ and $\g_2$ so that $\xi_2=e^{i\frac{\pi}{3}}\xi_1$. 
\end{itemize}

{\it We start with the case $a=1$, namely $\G_{\cal N}=\G_C$.} We
want to define a map 
$\varphi$ from $\C$ to $\Bn$ and then to extend $\varphi$ to
$\C\times\R$ by $\varphi(z,v)=\varphi(z)$, so that this extended map is
equivariant w.r.t. the action of $\G_C$. An obvious necessary condition is
that $\varphi:\C\fd\Bn$ must be 
equivariant w.r.t. the action of $T(\G_C)$ on $\C$. Another necessary condition is
that $\varphi$ should send $\C$ to the fixed points set ${\rm Fix}_0$ of
$\rho(\g_0)$ in $\Bn$. Indeed, for any $z\in\C$, $\g_0(z,0)$ belongs to
$\{z\}\times\R$ and $\varphi$ maps $\{z\}\times\R$ to the point $\varphi(z)$.
These two conditions are also clearly sufficient. 

So let $x_0\in\Bn$ be a fixed point of $\rho(\g_0)$ and set
$\varphi(0)=x_0$. 
Since $\g_i=(\xi_i,\nu_i)$, $i=1$ or $2$, commutes with $\g_0$,
the point $x_i=\rho(\g_i)x_0$ must also be fixed by $\rho(\g_0)$. 
Let $\sigma_{0i}$ be the geodesic arc in $\Bn$ joining $x_0$
to $x_i$. Note that 
${\rm Fix}_0$ is a convex subset of $\Bn$ and hence $\sigma_{0i}$ is included in
${\rm Fix}_0$. Let $\varphi$ map the segment $[0,\xi_i]$ onto $\sigma_{0i}$. We
then map the segment $[\xi_1,\xi_1+\xi_2]$ to $\rho(\g_1)\sigma_{02}$ and the
segment  $[\xi_2,\xi_1+\xi_2]$ to $\rho(\g_2)\sigma_{01}$. This is well defined
since
$\rho(\g_1)(x_2)=\rho(\g_1\g_2)(x_0)=
\rho(\g_2\g_1)\rho(\g_1^{-1}\g_2^{-1}\g_1\g_2)x_0=
\rho(\g_2\g_1)\rho(\g_0^{k})x_0=\rho(\g_2\g_1)x_0=\rho(\g_2)(x_1)$. Moreover,
because of the commutation of $\g_1$ and $\g_2$ with $\g_0$,
$\rho(\g_1)\sigma_{02}$ and $\rho(\g_2)\sigma_{01}$ are included in ${\rm
  Fix}_0$. 

Hence we get an equivariant map $\varphi$ from the boundary of a fundamental
domain of $T(\G_C)$ in $\C$ to ${\rm Fix}_0$ ($\varphi$ can be made smooth, for
example by taking it constant near $0$, $\xi_1$ and $\xi_2$). We can therefore
extend $\varphi$ to a $T(\G_C)$-equivariant map from $\C$ to ${\rm Fix}_0$.        

Define now $f:HS\times[0,+\infty)=\C\times\R\times[0,+\infty)\,\fd\,{\rm
Fix}_0\subset\Bn$ by $f(z,v,t)=\varphi(z)$. The map $f$ is $\rho$-equivariant and its
energy density is: 
$$
\|{\rm d}f\|^2_{(x,t)} 
= 
|{\rm d}f(\frac{\partial}{\partial t})|^2_{\varphi(x)}
+e^{2t}|{\rm d}f(\frac{\partial}{\partial v})|^2_{\varphi(x)}+
e^t\sum_{k=3}^{4}|{\rm d}f(e_k)|^2_{\varphi(x)}
= 
0+0+e^t\|{\rm d}\varphi\|^2_x
$$  
where $\|{\rm d}\varphi\|$ denotes the norm of the differential of
$\varphi:\C\fd {\rm Fix}_0$ computed with the metrics of
$\C\times\{0\}\subset HS\subset\Bm$ and ${\rm Fix}_0\subset\Bn$ induced
from $g_m$ and $g_n$.  

Therefore,
$
E_C(f)=\frac{1}{2}\int_0^{+\infty}\int_N \|{\rm d}\varphi\|^2e^{-t}dV_Ndt <
+\infty~. 
$
  
{\it Now, consider the cases where $a\neq 1$.} We want to proceed as we
just did, namely, we want to first construct a map $\varphi$ from $\C$ to $\Bn$
and then extend it to $HS$ by requiring that $\varphi(z,v)=\varphi(z)$. The two
conditions we mentioned are of course still necessary but we need to be more
careful because of the action of  $\G_C/\G_{\cal N}$. 

Let $\g_3$ be an element of $\G_C$ such that $\g_3\G_{\cal N}=a$. Then
$\g_3=(\xi_3,\nu_3,a)$ for some $\xi_3\in\C$ and $\nu_3\in\R$. It is easy
to check that $\g_0$, $\g_1$, $\g_2$ and $\g_3$ generate $\G_C$. 

The first thing to notice is that the point $\zeta=\frac{\xi_3}{1-a}$ is fixed by
the action of $\g_3$ on the $\C$-factor. Since we want $\varphi$
to be constant on $\{\zeta\}\times\R$, $\varphi$ must send
$\zeta$ to a fixed point of $\rho(\g_3)$. This can be done
because of the:        
 
\begin{claim}\label{ptfixecommun}
Let $\g=(\xi,\nu,b)\in\G_C$ be such that $b\neq 1$. Then
$\rho(\g)$ and $\rho(\g_0)$ have a common fixed point in $\Bn$.   
\end{claim}

\begin{demo} Since $\g$ and $\g_0$ commute, $\rho(\g)$ stabilizes
  the totally geodesic submanifold ${\rm Fix}_0$ of $\Bn$.  
Let $q$ be such that $b^q=1$. Computing, we get 
$
{\g}^q=(\xi,\nu,b)^q=((\sum_{k=0}^{q-1}b^k)\xi,v,b^q)=(0,v,1)
$ 
for some $v\in\R$.
Hence ${\g}^q$ belongs to $\G_{\cal Z}$: ${\g}^q$ is a
power of $\g_0$. The orbit under the group generated by $\rho(\g)$ of any
point in ${\rm Fix}_0$ must therefore be finite and this implies that the
action of $\rho(\g)$ on ${\rm Fix}_0$ has a fixed point. 
\end{demo}

With this in mind, it is now possible to complete the proof by constructing
$\varphi$ on the boundary $\partial F$ of a fundamental domain $F$ of the action
of $\G_C$ on the $\C$-factor. Since $\g_3$ acts on $\C$ by rotation around its
fixed point $\zeta$, we can choose a fundamental domain $G$ of the action of
$T(\G_{\cal N})$ on $\C$, centered at $\zeta$ and invariant by $\g_3$. For $F$
we then take a fundamental domain for the action of $\g_3$ on $G$.   

We do it in the case $a=e^{i\frac{\pi}{3}}$, the other
cases are handled similarly.

The lattice $T(\G_{\cal N})$ is generated by $\xi_1$ and
$\xi_2=e^{i\frac{\pi}{3}}\xi_1$. Let $G$ be the regular hexagon centered at
$\zeta$ with one vertex at the point $\zeta+\frac{1}{3}(\xi_1+\xi_2)$. $G$ is a
fundamental domain for the action of $T(\G_{\cal N})$ and it is invariant by
$\g_3$. Let then $F$ be the quadrilateral whose vertices are $\zeta$,
$\zeta+\frac{1}{2}\xi_1$, $\zeta+\frac{1}{3}(\xi_1+\xi_2)$ and
$\zeta+\frac{1}{2}\xi_2$. $F$ is clearly a fundamental
domain for the action of $\G_C$ on $\C$. See Fig.~1 for a picture. 

\begin{center}
\begin{picture}(0,0)%
\includegraphics{formfig3.pstex}%
\end{picture}%
\setlength{\unitlength}{1579sp}%
\begingroup\makeatletter\ifx\SetFigFont\undefined%
\gdef\SetFigFont#1#2#3#4#5{%
  \reset@font\fontsize{#1}{#2pt}%
  \fontfamily{#3}\fontseries{#4}\fontshape{#5}%
  \selectfont}%
\fi\endgroup%
\begin{picture}(13584,9830)(-454,-10083)
\put(7588,-3916){\makebox(0,0)[lb]{\smash{\SetFigFont{8}{9.6}{\rmdefault}{\mddefault}{\updefault}
$\xi_2$
}}}
\put(8244,-5213){\makebox(0,0)[lb]{\smash{\SetFigFont{8}{9.6}{\rmdefault}{\mddefault}{\updefault}
$\xi_1$
}}}
\put(6516,-5697){\makebox(0,0)[lb]{\smash{\SetFigFont{8}{9.6}{\rmdefault}{\mddefault}{\updefault}
$\zeta$
}}}
\put(7070,-5109){\makebox(0,0)[lb]{\smash{\SetFigFont{8}{9.6}{\rmdefault}{\mddefault}{\updefault}
F
}}}
\end{picture}
 \\  
Fig.~1  
\end{center}

Let now $x_0\in\Bn$ be a fixed point of both $\g_0$ and $\g_3$ (such a point
exists by Claim~\ref{ptfixecommun}). Set $\varphi(\zeta)=x_0$. 

The point $\zeta+\frac{1}{2}\xi_1$ is fixed by $\g_1\g_3^3$, hence it must be sent
by $\varphi$ to a fixed point of $\rho(\g_1\g_3^3)$. It follows from Claim~\ref{ptfixecommun} that $\rho(\g_1\g_3^3)$ and $\rho(\g_0)$ have a common fixed
point, say $x_1$. Let $\varphi$ send the vertex $\zeta+\frac{1}{2}\xi_1$ to
$x_1$ and the edge $[\zeta,\zeta+\frac{1}{2}\xi_1]$ of $F$ to the geodesic arc
$\sigma_{01}$ joining $x_0$ to $x_1$ in ${\rm Fix}_0$. Similarly, the vertex
$\zeta+\frac{1}{3}(\xi_1+\xi_2)$ is a fixed point of $\g_2\g_3^4$ and we let
$\varphi$ map it to a fixed point $x_2$ of $\rho(\g_2\g_3^4)$ in ${\rm
  Fix}_0$. We map the edge
$[\zeta+\frac{1}{2}\xi_1,\zeta+\frac{1}{3}(\xi_1+\xi_2)]$ to the geodesic arc
$\sigma_{12}$ joining $x_1$ and $x_2$ in ${\rm Fix}_0$. 

Now the edge $[\zeta,\zeta+\frac{1}{2}\xi_2]$ is the image of
$[\zeta,\zeta+\frac{1}{2}\xi_1]$ under $\g_3$ so we must map it to
$\rho(\g_3)(\sigma_{01})$. In the same way, $[\zeta+\frac{1}{2}\xi_2,
\zeta+\frac{1}{3}(\xi_1+\xi_2)]$ is the image of $[\zeta+\frac{1}{2}\xi_1,
\zeta+\frac{1}{3}(\xi_1+\xi_2)]$ by $\g_2\g_3^4$ and we must therefore map it to
$\rho(\g_2\g_3^4)(\sigma_{12})$. 
These definitions of $\varphi$ agree at the point $\frac{1}{2}\xi_2$. Indeed,
a simple computation shows that there exists $q$ such that
$\g_2\g_3=\g_0^q\g_3\g_1$ and therefore,  
$\rho(\g_2\g_3^4)x_1=\rho(\g_2\g_3\g_3^3)x_1=\rho(\g_0^q\g_3\g_1\g_3^3)x_1=
\rho(\g_0^q\g_3)\rho(\g_1\g_3^3)x_1=\rho(\g_3)\rho(\g_0^q)x_1=\rho(\g_3)x_1$.
Hence $\varphi$ is well defined on $\partial F$. By construction, $\varphi$
is equivariant w.r.t. $\g_0$ and the face-pairings $\g_3$ and $\g_2\g_3^4$
which generate the whole group $\G_C$. 

The construction of $\varphi$ and $f$ then goes on as in the case $a=1$.

\vspace{0.2cm}

In this way we obtain a section $f_i$ of the bundle $M\times_\rho\Bn$ on
each cusp $C_i$ of $M$. This section can be extended to a section $f$
defined on the whole manifold $M$ and since the energy of
$f_i:C_i\fd M\times_\rho\Bn$ is finite for each $i$, the energy of
$f:M\fd M\times_\rho\Bn$ is finite and we are done.     
\end{demo}

\section{Pluriharmonicity and consequences}\label{pluri}

In this section, we study the properties of finite energy harmonic maps
$\Bm\fd\Bn$ which are equivariant w.r.t. a representation $\rho$ of a
torsion-free lattice $\G<{\rm PU}(m,1)$ into ${\rm PU}(n,1)$. 

\subsection{Pluriharmonicity}\hfill 

\begin{theo}\label{pluriharmonic}
  Let $f:\Bm\fd\Bn$ be a $\rho$-equivariant harmonic map of finite
  energy.   
  Then $f$ is pluriharmonic, namely, the $J$-invariant part ${(\nab\df)}^{1,1}$
  of $\nab\df$ vanishes identically. Moreover the complexified sectional
  curvature of $\Bn$ is zero on $\df(T^{1,0}\Bm)$.
\end{theo}

We first prove a general Bochner-type formula due to Mok, Siu and Yeung
(cf. \cite{MSY93}) in case $\G$ is a uniform lattice. We state it in the case of
maps $\Bm\fd\Bn$ but it is valid in the more general setting of
equivariant maps from an irreducible (rank 1) symmetric space of non-compact type
to a negatively curved manifold, as can be seen from the proof. Our
exposition follows~\cite{Pa95}.


Let $R^m$ and $R^n$ be the curvature tensors of $g_m$ and $g_n$, and $Q$ be
any parallel tensor of curvature type on $\Bm$. For $h$ a symmetric
2-tensor with values in a vector bundle over $\Bm$, define $(\Qrond h)(X,Y) =
\mbox{tr}(W\longmapsto h(Q(W,X)Y,W))$.   

Remark that if $f$ is a $\rho$-equivariant map $\Bm\fd\Bn$, then, since $Q$ is
parallel and $\rho(\Gamma)$ acts by isometries on $\Bn$, the $\R$-valued
functions $\la Q,f^\star R^n\ra$ and $\la\Qrond\nab\df,\nab\df\ra$ on $\Bm$ 
are in fact $\rho$-invariant and hence can be considered as functions on
$M=\Gamma\backslash\Bm$: 

\begin{prop}\label{bochnerformula}
Let $f$ be a $\rho$-equivariant harmonic map of finite energy from $\Bm$ to
$\Bn$ and $Q$ a parallel tensor of curvature type on $\Bm$. Then,
\begin{equation}
\int_M\left[\la\Qrond\nab\df,\nab\df\ra -\frac{1}{2}\la Q,f^\star R^n\ra\right]dV_m = 
-\frac{1}{4m}\int_M \la Q,R^m\ra\|\df\|^2 \,dV_m~,\tag{\mbox{$\diamondsuit$}}
\end{equation}
where if $\G$ is non-uniform, that is if $M$ is non-compact,  the left-hand side
should read $\lim_{R\pfd\infty}\int_M \eta_R[\la\Qrond\nab\df,\nab\df\ra 
-\frac{1}{2}\la Q,f^\star R^n\ra]dV_m$, for $\{\eta_R\}$ a well-chosen
family of cut-off functions on $M$.  
\end{prop}

\begin{demo} Let us first assume that $M$ is compact. 
All the computations will be made in a normal coordinates system. 

By definition, $(\Qrond\nab\df)(X,Y)=\sum_k
(\nab\df)(Q(e_k,X)Y,e_k)$. Since $\df$ is closed, i.e. $\nab\df$
is symmetric, and $Q$ is parallel, we have in fact 
$$
\Qrond\nab\df(X,Y)=\sum_k
(\nab_{e_k}\df)(Q(e_k,X)Y)= 
\sum_k (\nab_{e_k}\df\circ Q)(e_k,X)Y 
= -\nab^\star (\df\circ Q)(X,Y)~,
$$
where $\nab^\star $ is the formal adjoint of $\nab$: if $T$ is a
$(p+1)$-tensor,
$(\nab^\star T)(X_1,\ldots,X_p):=-\mbox{tr}(W\longmapsto
(\nab_WT)(W,X_1,\ldots,X_p))$.

Integrating this relation over $M$ (we assumed $M$ compact), we get  
$$
\int_M \la\Qrond\nab\df,\nab\df\ra dV_m 
= 
-\int_M\la\nab^\star (\df\circ Q),\nab\df\ra dV_m 
= 
-\int_M\la\df\circ Q,\nab^2\df\ra dV_m 
$$
where $\nab^2\df$ is the 3-tensor $\nab(\nab\df)$.

Using that $Q$, and hence $\df\circ Q$, is skew-symmetric in its first two
variables, one checks that 
$$
\la\df\circ Q,\nab^2\df\ra 
=-\frac{1}{2} 
[\la\df\circ Q,f^\star R^n\ra-\la\df\circ Q,\df\circ R^m\ra]~.
$$
We have $\la\df\circ Q,f^\star R^n\ra=\la Q,f^\star R^n\ra$, 
where in the  r.h.s. $f^\star R^n$ and $Q$ are considered as (4,0)-tensors.  
Moreover, computations show that
$$
\la \df\circ Q,\df\circ R^m\ra 
=\sum_{a,b}\frac{1}{2}\Big(\la\iota_{e_a}Q,\iota_{e_b}R^m\ra
+\la\iota_{e_b}Q,\iota_{e_a}R^m\ra\Big)f^\star\gn(e_a,e_b)~,
$$
where $\iota$ denotes interior product. Now, since $M$ is locally symmetric, the
symmetric 2-tensor $\t$ given by   
$$\t(X,Y)=\frac{1}{2}(\la\iota_{X}Q,\iota_{Y}R^m\ra
+\la\iota_{Y}Q,\iota_{X}R^m\ra)
$$
is parallel. Thus it must be proportionnal to $\gm$ ($M$ is locally
irreducible): $\t=\frac{1}{2m}({\rm tr}_{\gm}\t)\,\gm$. Now, ${\rm
  tr}_{\gm}\t=\la Q,R^m\ra$ and $\la \gm,f^\star\gn\ra =\|\df\|^2$, so
that
$
\la \df\circ Q,\df\circ R^m\ra = \frac{1}{2m}\la Q,R^m\ra \|\df\|^2
$ 
and hence
$$
\int_M\la \Qrond\nab\df,\nab\df\ra dV_m = 
\frac{1}{2}\int_M\left[\la Q,f^\star R^n\ra-\frac{1}{2m}\la Q,R^m\ra\|\df\|^2\right]dV_m~.
$$
This ends the proof in the compact case. 

\vspace{0.2cm}

Now assume $M$ is non-compact of finite volume. The only global step in
the preceding proof is the initial integration by parts. Thus we only have
to show that this can be done in the finite volume case. 
We mimic the argument given by Corlette in~\cite{Corlette92}.

As mentionned earlier, $M$ is the union of a compact manifold with boundary
$M_0$ and of a finite 
number of pairwise disjoint cusps $C_i$, each diffeomorphic to a compact
$(2m-1)$-manifold $N_i$ times $[0,+\infty)$. For each $i$, let $t_i$ be the
parameter in the $[0,+\infty)$ factor. 

For $R>1$, we define a cut-off function $\eta_R$ on $M$ in the following
manner. Take a smooth function $\eta$ on $[0,+\infty)$ identically equal to
1 on $[0,1]$ and to 0 on $[2,+\infty)$. Set 
$$
\eta_R(x)\,= 
\left\{
  \begin{array}{ll}
1 & \mbox{if } x\in M_0, \\
\ds \eta\Big(\frac{t_i}{R}\Big) & \mbox{if } x\in C_i.  
  \end{array}
\right.
$$ 
Since $\eta_R$ is a horofunction along each cusp (cf. \cite{Helgason84},
II.3.8), the absolute value $|\Delta \eta_R|$ of its Laplacian is bounded
independently of $R$. Moreover, the norm $\|{\rm d}\eta_R\|$ of its differential
is bounded by a constant times $\frac{1}{R}$.      

Introducing $\eta_R$ in the integration by parts, we obtain 
$$
\begin{array}{rcl}
\displaystyle \int_M \eta_R\,\la\Qrond\nab\df,\nab\df\ra dV_m 
           & = & \ds -\int_M \la\nab^\star (\df\circ Q),\eta_R\,\nab\df\ra dV_m \\
           & = & \ds -\int_M \la\df\circ Q,\eta_R\,\nab^2\df + {\rm
             d}\eta_R\otimes\nab\df \ra dV_m \\
           & = & \ds -\int_M \eta_R\, \la\df\circ Q,\nab^2\df\ra dV_m -\int_M
           \la\df\circ Q, {\rm d}\eta_R\otimes\nab\df \ra dV_m \\
           & = & \ds \frac{1}{2}\int_M \eta_R\,\left[\la Q,f^\star R^n\ra
           -\frac{1}{2m}\,\la Q,R^m\ra\|\df\|^2\right] dV_m\\ 
           &   & \hfill \ds-\int_M\la\df\circ Q, {\rm d}\eta_R\otimes\nab\df \ra dV_m~. 
\end{array}
$$
Thus,
$$
\begin{array}{rcl}
\ds \int_M \eta_R\left[\la\Qrond\nab\df,\nab\df\ra-\frac{1}{2}\la Q,f^\star
  R^n\ra\right] dV_m
& = & \ds  -\frac{1}{4m}\int_M \eta_R\la Q,R^m\ra\|\df\|^2 dV_m \\
&   & \hfill \ds -\int_M\la\df\circ Q, {\rm d}\eta_R\otimes\nab\df \ra dV_m~. 
\end{array}
$$
The tensors $Q$ and $R^m$ are parallel and hence $\la Q,R^m\ra$ is constant on
$M$. Therefore, $\la Q,R^m\ra\|\df\|^2$ is integrable and the first term in the
r.h.s. goes to $-\frac{1}{4m}\int_M \la Q,R^m\ra\|\df\|^2 dV_m$ as $R$ goes to
infinity. On the other hand we have $\|{\rm d}\eta_R\|\leq \frac{C}{R}$ for some
constant $C$ independent of $R$ and hence
$$
\begin{array}{rcl}
\ds {\Big(\int_M \la\df\circ Q, {\rm d}\eta_R\otimes\nab\df \ra dV_m\Big)}^2
& \leq & \ds
\Big(\int_M \|\df\circ Q\|^2 dV_m \Big)\,\Big(\int_M\|{\rm
  d}\eta_R\otimes\nab\df\|^2\,dV_m\Big)\\
& \leq & \ds 
\Big(\int_M \frac{1}{2m}\,\|Q\|^2\|\df\|^2 dV_m \Big)\,\Big(\int_M 2m\,\|{\rm
  d}\eta_R\|^2\,\|\nab\df\|^2\,dV_m\Big)\\
& \leq & \ds 
\frac{C^2}{R^2}\, \Big(\int_M \|Q\|^2\|\df\|^2 dV_m \Big)\,\Big(\int_M
\|\nab\df\|^2\,dV_m\Big)~. 
\end{array}
$$
Since $Q$ is parallel, $\|Q\|$ is constant and $\int_M \|Q\|^2\|\df\|^2 dV_m$
is finite. 

The next lemma implies that $\lim_{R\pfd\infty}\int_M\la\df\circ
Q,{\rm d}\eta_R\otimes\nab\df\ra dV_m = 0$ and therefore ends the proof of
Proposition~\ref{bochnerformula}. 
\end{demo}  

\begin{lemma}\label{nabladfl2}
$ \|\nab\df\|$ belongs to $L^2(M)$: $\int_M \|\nab\df\|^2\,dV_m<+\infty$.     
\end{lemma}

\begin{demo}
Because the energy density $e$ of $f$ is integrable on
$M$, and using Green's formula, we see that  
$$
\int_M (\Delta e)\,\eta_R\,dV_m = \int_M e\,(\Delta \eta_R)\,dV_m 
$$
is bounded independently of $R$. Now, since we assumed that $f$ is harmonic, the
Bochner-type formula of Eells-Sampson (\cite{ES64}) reads: 
$$
\Delta e = -\|\nab\df\|^2+{\rm Scal}(f^\star R^n)-\la \df\circ {\rm Ric}^m,\df\ra~,
$$
where ${\rm Scal}(f^\star R^n)$ denotes the scalar contraction of
the curvature tensor $f^\star R^n$ and ${\rm Ric}^m$ is the Ricci tensor of $g_m$
seen as an endomorphism of $T\Bm$. Since $\Bn$ is negatively curved and
${\rm Ric}^m =-\frac{1}{2}(m+1){\rm Id}$, we get $\|\nab\df\|^2 \leq -\Delta e +
(m+1)\,e$ and thus $\|\nab\df\|^2$ is integrable. 
\end{demo}

Let us call $I$, resp. $I_\C$, the (3,1)-tensor of curvature type on
$M=\G\backslash\Bm$ (or on $\Bm$) given by $I(X,Y)Z=\gm(X,Z)Y-\gm(Y,Z)X$,
resp. $I_\C(X,Y)Z=\frac{1}{4}(I(X,Y)Z+I(JX,JY)Z+2\gm(JX,Y)JZ)$, for all
$X,Y,Z\in T\Bm$. The curvature tensor $R^m$ of $M$ (or of $\Bm$) is just
$-I_\C$. Both $I$ and 
$I_\C$ are parallel tensors, and in fact they form a basis of the space of
parallel tensors of curvature type on $M$. $I$ and $I_\C$ will also denote the
corresponding (4,0)-tensors. 

We will apply the Bochner-type formula $(\diamondsuit)$ to the parallel
tensor of curvature type $Q=I_\C-I$. 

\begin{lemma}\label{Jinv}
  Let $f:\Bm\fd\Bn$ be a harmonic map and let $Q=I_\C-I$.\\ 
Then $\la\Qrond\nab\df,\nab\df\ra=-\frac{3}{2}\|(\nab\df)^{(1,1)}\|^2$,  where
$(\nab\df)^{(1,1)}$ is the $J$-invariant part of $\nab\df$: for all $X,Y\in T\Bm$,
$(\nab\df)^{(1,1)}(X,Y):=\frac{1}{2}[\nab\df(X,Y)+\nab\df(JX,JY)]$.  
\end{lemma}

\begin{demo}
A straightforward computation shows that for $h$ a symmetric 2-tensor
taking values in $f^\star T\Bn$,
$\stackrel{\,\circ}{I_{~}}\hspace{-2pt}h=h-g_m\,{\rm tr}_{g_m}h$ and 
$\stackrel{\,\circ}{I_\C}h(X,Y)=\frac{1}{4}[h(X,Y)-3h(JX,JY)-g_m(X,Y)\,{\rm
  tr}_{g_m}h]$. Therefore, since ${\rm tr}_{g_m}\nab\df=0$,
$\Qrond\nab\df=-\frac{3}{2}(\nab\df)^{(1,1)}$. The
decomposition of a 2-tensor in $J$-invariant and $J$-skew-invariant parts is
orthogonal, hence the result.     
\end{demo}

\begin{lemma}\label{ortho}
  $I_\C$ is the orthogonal projection of $I$ onto the space of Kähler
  curvature type tensors, namely, the space of tensors of curvature type $T$
  such that $T(X,Y)JZ=J(T(X,Y)Z)$, for all $X,Y,Z\in T\Bm$.
\end{lemma}

\begin{demo} Since $I_\C$ is clearly of Kähler curvature type, it remains to show
  that $I_\C-I$ is orthogonal to all tensors of Kähler curvature type. 
  Simple computations shows that if $T$ is any tensor of curvature
  type, $\la I,T \ra = 2\,{\rm Scal}(T)$, whereas  
$$
\la I_\C,T \ra = \frac{1}{2}{\rm Scal}(T)-
\frac{1}{2}\sum_{k,l=1}^{2m}\Big(T(e_k,Je_l,Je_k,e_l)+T(e_k,Je_k,Je_l,e_l)\Big),
$$ 
for $\{e_k\}$ an orthonormal basis of $TM$. It is then easy to check that if $T$
is moreover of Kähler type, this last formula reduces to $\la I_\C,T \ra  =
2\,{\rm Scal}(T)$, hence the result. 
\end{demo}

Let us recall what the {\it complexified sectional curvature} of 
a Hermitian manifold $(N,g,J)$ is:
if $E$ and $F$ are two vectors of the complexified tangent space
$T^\C N=TN\otimes_\R\C$ of $N$ then the complexified sectional curvature of the
2-plane they span is defined to be $R^N(E,F,\ov E,\ov F)$ where $R^N$ is the
curvature tensor of $g$ extended by $\C$-linearity to $T^\C N$. Despite its
name, the complexified sectional curvature takes real values.

If $T$ is a tensor of curvature type, we define its {\it complexified scalar
curvature} ${\rm Scal}_\C(T)$ as follows: ${\rm Scal}_\C(T):=\sum_{k,l=1}^m
T(\zeta_k,\zeta_l,\ov \zeta_k,\ov \zeta_l)$, for $\{\zeta_k\}$ an orthonormal
basis of the $(1,0)$-part of $T^\C \Bm$.

Using the formulae given in the proof of the previous lemma, one gets
\begin{lemma}\label{complexified}
  $\la I_\C-I,T\ra= -6\, {\rm Scal}_\C(T)$.
\end{lemma}

We are now ready to prove Theorem~\ref{pluriharmonic}. Recall that
$Q=I_\C-I$. First,
Lemma~\ref{ortho} implies that the right-hand side in the Bochner-type formula
($\diamondsuit$) vanishes. Next, it follows from Lemma \ref{Jinv} and Lemma
\ref{complexified} that   
$$
\int_M\eta_R\left[\la\Qrond\nab\df,\nab\df\ra-\frac{1}{2}\la Q,f^\star
  R^n\ra\right]dV_m= -\frac{3}{2}\int_M\eta_R\left[\|(\nab\df)^{(1,1)}\|^2-2\,{\rm
    Scal}_\C(f^\star R^n)\right]dV_m 
$$ 
for any $R>1$. Thus, formula ($\diamondsuit$) reads:
$$
\lim_{R\pfd\infty}\int_M\eta_R\left[\|(\nab\df)^{(1,1)}\|^2-2\,{\rm
    Scal}_\C(f^\star R^n)\right]dV_m\,=\,0~.
$$

It is known that, since the sectional curvature of $(\Bn,g_n)$ is pinched
between $-1$ and $-\frac{1}{4}$, its complexified sectional curvature is
non-positive (see for example~\cite{Her91}). Therefore, ${\rm
  Scal}_\C(f^\star R^n)$ being a mean of 
complexified sectional curvatures of $\Bn$, it is non-positive. Thus $R\longmapsto 
\int_M\eta_R\left[\|(\nab\df)^{(1,1)}\|^2-2\,{\rm Scal}_\C(f^\star
  R^n)\right]dV_m$ is a non-negative non-decreasing function whose limit as $R$
goes to infinity is zero. It follows at once that $(\nab\df)^{(1,1)}$ vanishes
identically, that is, $f$ is pluriharmonic. Finally, we also have ${\rm
  Scal}_\C(f^\star R^n)=0$ everywhere and this implies that $R^n(f_\star
\zeta_k,f_\star \zeta_l,\ov{f_\star \zeta_k},\ov{f_\star \zeta_l})=0$ for all
$k,l=1,\ldots,n$. 
Theorem \ref{pluriharmonic} is proved.

\subsection{Holomorphicity of ``high'' rank harmonic maps}\hfill

Let $f:\Bm\fd\Bn$ be a finite energy harmonic map, equivariant
w.r.t a representation 
of a torsion-free lattice $\G<{\rm PU}(m,1)$ into ${\rm PU}(n,1)$. In this
section, we exploit 
the full strength of Theorem~\ref{pluriharmonic} to prove  
a result that was first obtained by J.~A.~Carlson and D.~Toledo in~\cite{CT89}
in the case $\G$ is cocompact and $\rho(\G)$ is discrete in ${\rm PU}(n,1)$ or
in a more general target Lie group $G$. Their
proof relies on a careful study of maximal abelian subalgebras of the
complexification of the Lie algebra of $G$. In our setting, the simple
form of the curvature tensor of $(\Bn,g_n)$ allows a more elementary proof
that we give for completeness.

\begin{prop}\label{holom}
  Let $f:\Bm\fd\Bn$ be a finite energy harmonic map equivariant
  w.r.t. a
  representation $\rho$ of the torsion-free lattice $\G<{\rm PU}(m,1)$ in
  ${\rm PU}(n,1)$. If the real rank of $f$ is at least 3 at some point,
  then $f$ is holomorphic or anti-holomorphic.  
\end{prop}

Before proving this proposition, we introduce some notations that will be
needed in the proof and later on in the paper.

For $l=m,n$, let
$T^\C\Bl=T\Bl\otimes_\R\C$ be the complexification of $T\Bl$ and $T^\C
\Bl=T^{1,0}\Bl\oplus T^{0,1}\Bl$ be its decomposition in (1,0) and (0,1)
part. We extend the differential of $f$ by $\C$-linearity and still write 
$\df:T^\C \Bm\fd T^\C \Bn$ (if the distinction is necessary we will use
${\rm d}^\C f$). Its components are 
$$ 
\partial^{1,0}f:T^\C \Bm\fd T^{1,0} \Bn~,
$$
$$
\partial^{0,1}f:T^\C \Bm\fd T^{0,1} \Bn~.
$$

We extend $g_n$ by $\C$-linearity to $T^\C\Bn$. We will sometimes write
$(X,Y)=g_n(X,Y)$ and $|X|^2=(X,\ov X)$ for $X,Y\in T^\C\Bn$.  

If $(e_i)_{1\leq i\leq 2m}=(e_\alpha,Je_\alpha)_{1\leq\alpha\leq m}$ is an
orthonormal $\R$-basis of $T\Bm$, we set
$z_\alpha=\frac{1}{2}(e_\alpha-iJe_\alpha)$. 
$(z_\alpha)_{1\leq\alpha\leq m}$ is an orthogonal $\C$-basis of $T^{1,0}\Bm$. 

To lighten the notations, we will sometimes use $f_\a$ instead of 
$\partial^{1,0}f(z_\alpha)$ and $f_{\ov\b}$ instead of $\partial^{1,0}f(\bar
z_\beta)$, so that
$\partial^{0,1}f(\bar z_\alpha)=\bar f_\a$ and
$\partial^{0,1}f(z_\beta)=\bar f_{\ov\b}$. 

In the sequel, we will often restrict $\partial^{1,0}f$, resp. $\partial^{0,1}f$,
to $T^{1,0}\Bm$ and consider them as sections of ${\rm Hom}(T^{1,0}\Bm,f^\star
T^{1,0}\Bn)$, resp. ${\rm Hom}(T^{1,0}\Bm,f^\star T^{0,1}\Bn)$. 
We will call $e'(f)$, resp. $e''(f)$, the square of the norm of $\partial^{1,0}f$,
resp. $\partial^{0,1}f$, namely:

$$
\begin{array}{rcl}
e'(f)& := & \ds \|\partial^{1,0}f\|^2=2 \sum_{\a=1}^m
\gn(\partial^{1,0}f(z_\a),\ov{\partial^{1,0}f(z_\a)})=2\sum_{\a=1}^{m}|f_\a|^2~, \\ 
e''(f) & := & \ds \|\partial^{0,1}f\|^2= 2\sum_{\a=1}^m
\gn(\partial^{0,1}f(z_\a),\ov{\partial^{0,1}f(z_\a)})=2\sum_{\a=1}^{m}|f_{\ov\a}|^2~. 
\end{array}
$$

Note that with these definitions, the energy density of $f$ is given by
$e(f):=\frac{1}{2}\|\df\|^2=e'(f)+e''(f)$. Again, we will often abreviate $e'(f)$ and
$e''(f)$ to $e'$ and $e''$ when no confusion is possible.   

\medskip


\begin{demo}[of Proposition~\ref{holom}]
   Theorem \ref{pluriharmonic} shows that $f^\star
   R^n(\za,\zb,\zab,\zbb)=0$ for all $\a,\b\in\{1,\ldots,m\}$. Let us see
   what this implies in more details.


$$\begin{matrix}
f^\star R^n(z_\alpha,z_\beta,\bar z_\alpha,\bar z_\beta)& = & R^n(f_\a+\bar f_{\ov\a},f_\b+\bar f_{\ov\b},\bar f_\a+f_{\ov\a},\bar f_\b+f_{\ov\b}) \hfill\cr
& = & -\frac{1}{2}\ds
\Big[|f_\a|^2|f_{\ov\b}|^2+|f_{\ov\a}|^2|f_\b|^2-(
f_{\ov\a},\bar f_{\ov\b})( f_\a,\bar f_\b)\hfill\cr 
&&\hfill -(f_\b,\bar f_\a)( f_{\ov\b},\bar f_{\ov\a})-(
f_\a,\bar f_{\ov\b})( f_{\ov\a},\bar f_\b)-( f_\b,\bar f_{\ov\a})(
f_{\ov\b},\bar f_\a)\hfill\cr 
&& \hfill+( f_\b,\bar f_{\ov\a})( f_{\ov\a},\bar f_\b)+( f_\a,\bar
f_{\ov\b})( f_{\ov\b},\bar f_\a)\Big]\cr
&=& -\frac{1}{2}\big(\|f_\a\wedge \bar f_{\ov\b}-f_\b\wedge\bar f_{\ov\a}\|^2+|(
f_\a,\bar f_{\ov\b})-( f_\b,\bar f_{\ov\a})|^2\big)\hfill\cr
\end{matrix}$$
Therefore, for all $\alpha$ and $\beta$, we have $f_\a\wedge \bar
f_{\ov\b}=f_\b\wedge\bar f_{\ov\a}$. 

Suppose that the families $(f_\a)_{1\leq\alpha\leq m}$ and
$(f_{\ov\b})_{1\leq\beta\leq m}$ are both of rank less than or equal to
1. We may assume that for all $\a$ there exits $\lambda_\a$ such that
$f_\a=\lambda_\a f_1$ and that for some $k$ and for all $\b$ there exists
$\mu_\b$ such that $f_{\ov \b}=\mu_\b f_{\ov k}$. 

If $f_1=0$, then for all $\a$, ${\rm d}^\C f(\za)=\bar f_{\ov\a}$ and
${\rm d}^\C f(\zab)=f_{\ov\a}$. Therefore the complex rank of ${\rm
  d}^\C f$ is at most 2, namely the real rank of ${\rm d}f$ is at most
2. The same is true if $f_{\ov k}=0$.

If both $f_1$ and  $f_{\ov k}$ are non zero, then from the fact that
$f_1\wedge\bar f_{\ov\a}=f_\a\wedge\bar f_{\ov 1}$ we deduce that
$\ov{\mu_\a}f_1\wedge\bar f_{\ov k}=\lambda_\a \ov{\mu_1}f_1\wedge\bar
f_{\ov k}$, that is, $\ov{\mu_\a}=\lambda_\a \ov{\mu_1}$. Then,  
$$
{\rm  d}^\C f(\za)=f_\a+\bar f_{\ov\a}=\lambda_\a f_1+\ov{\mu_\a}\bar f_{\ov k}=
\lambda_\a f_1+\lambda_\a \ov{\mu_1}\bar f_{\ov k}= \lambda_\a(f_1+ \bar
f_{\ov 1})=\lambda_\a\,{\rm d}^\C f(z_1)~.
$$ 
Hence the family $({\rm d}^\C f(\za))_{1\leq\alpha\leq m}$ has rank $\leq
1$. This also holds for the family $({\rm d}^\C f(\zbb))_{1\leq\b\leq m}$
and we conclude that the real rank of $f$ is again less than or
equal to 2. 

In any case, we see that if the real rank of $f$ is at least 3 at some
point, then the rank of one of the families $(f_\a)_{1\leq\alpha\leq m}$ and
$(f_{\ov\b})_{1\leq\beta\leq m}$ is at least 2 at this point.

Suppose now that at some point of $\Bm$, the rank of the family
$(f_\a)_{1\leq\alpha\leq m}$ is at least 2, for example $f_1\wedge
f_2\not=0$. Then $f_1\wedge \bar f_{\ov 2}=f_2\wedge\bar f_{\ov 1}$ implies
$f_{\ov 1}=f_{\ov 2}=0$. From $f_1\wedge \bar f_{\ov
  \gamma}=f_\gamma\wedge\bar f_{\ov 1}$,
we conclude that $f_{\ov\gamma}=0$ for all $1\leq\gamma\leq m$,
i.e. $e''=0$ at this point.  

In the same way, if the rank of the family $(f_{\ov\b})_{1\leq\beta\leq m}$
is at least 2, then $e'=0$. 

Finally, since $f$ is pluriharmonic and the complexified sectional
curvature of $\Bn$ is zero on $\df(T^{1,0}\Bm)$ (see Theorem \ref{pluriharmonic} above),
   and because $\Bn$ is a Kähler symmetric space, it is known that
   $\partial^{1,0}f$, resp. $\partial^{0,1}f$, are holomorphic sections of
   the holomorphic bundles 
   ${\rm Hom}(T^{1,0}\Bm,f^\star T^{1,0}\Bn)$, resp. ${\rm
     Hom}(T^{1,0}\Bm,f^\star T^{0,1}\Bn)$ (\cite{CT89}, Theorem~2.3). So they have a
   generic rank on $\Bm$. Therefore, if for example the family
   $(f_\a)_{1\leq\alpha\leq m}$ has rank at least 2 at some point, it has
   rank at least 2 on a dense open subset of $\Bm$ and so $e''=0$ on a
   dense open subset of $\Bm$, hence everywhere, and $f$ is
   holomorphic. Similarly, if ${\rm rk}(f_{\ov\b},\,{1\leq\beta\leq m})\geq
   2$ at some point, $f$ is antiholomorphic. 
\end{demo}

\subsection{Some technical lemmas}\label{lemmestechniques}\hfill

If $f:\Bm\fd\Bn$ is a $\rho$-equivariant pluriharmonic map whose rank is at
most 2 everywhere, $f$ needs not be holomorphic nor antiholomorphic. 
Nevertheless, pluriharmonicity has other consequences that will be useful
later. 
Namely, if $e'(f)$ and $e''(f)$ are the previously defined squared norms
of $\duzf$ and $\dzuf$, we have 
$$
\la f^\star\omega_n,\omega_m\ra:=
\sum_{i,j=1}^{2m}f^\star\omega_n(e_i,e_j)\omega_m(e_i,e_j)=2(e'(f)-e''(f))~,
$$  
as it is easy to check. Because of this, some results on the energies $e'(f)$ and
$e''(f)$ will be needed in the proof of Theorem~\ref{burgio} and we shall prove
them in this section. 

The results stated here were obtained in complex dimension~1 in~\cite{To79}
and~\cite{Wo79}. The proofs in the general case (see also \cite{Li70}) go
almost exactly as in the case $m=1$ and they are given only for
completeness and to fix the notations. 

We will work on the complexifications of the tangent spaces of $\Bm$ and
$\Bn$ and therefore we extend all needed sections, tensors and operators defined
on the real tangent spaces by $\C$-linearity to these complexifications.

Since $\duzf$ can be considered as a section of ${\rm Hom}(T^\C\Bm,f^\star
T^\C\Bn)$, we can define its covariant derivative $\nab\duzf\in{\rm
  Hom}(T^\C\Bm\otimes T^\C\Bm,f^\star T^\C\Bn)$.
It follows easily from the fact that $(\Bn,g_n)$ is Kähler
that $\nab\duzf$ belongs in fact to ${\rm
  Hom}(T^\C\Bm\otimes T^\C\Bm,f^\star T^{1,0}\Bn)$. We will call $\nab'\duzf\in{\rm
  Hom}(T^{1,0}\Bm\otimes T^{1,0}\Bm,f^\star T^{1,0}\Bn)$ its restriction to
$T^{1,0}\Bm\otimes T^{1,0}\Bm$.
We define  $\nab'\dzuf\in{\rm
  Hom}(T^{1,0}\Bm\otimes T^{1,0}\Bm,f^\star T^{0,1}\Bn)$ and
$\nab'\df=\nab'\duzf+\nab'\dzuf\in{\rm 
  Hom}(T^{1,0}\Bm\otimes T^{1,0}\Bm,f^\star T^\C\Bn)$ similarly.
Note that  $\|\nab'\df\|$, $\|\nab'\duzf\|$ and $\|\nab'\dzuf\|$ belong to
$L^2(M)$ because $\|\nab\df\|$ does (Lemma \ref{nabladfl2}).  

In the entire section, $f$ denotes a $\rho$-equivariant {\it pluriharmonic} map $\Bm\fd\Bn$.

\begin{lemma}\label{dep} We have
$$\frac{1}{4}\Delta e'=-\frac{1}{2}\|\nabla'\partial^{1,0}f\|^2-2R'+\frac{m+1}{4}\,e'
\;\;\mbox{ and }\;\;
\frac{1}{4}\Delta e''=-\frac{1}{2}\|\nabla'\partial^{0,1}f\|^2-2R''+\frac{m+1}{4}\,e''$$
where
$$
\begin{array}{rcl}
R'&=& \sum_{\a,\b} R^n(\df(\za),\df(\zab),\duzf(\zb),\ov{\duzfb})~,\\
R''&=&\sum_{\a,\b} R^n(\df(\za),\df(\zab),\dzuf(\zb),\ov{\dzufb})~.
\end{array}
$$
\end{lemma}

\begin{demo}
We make the computation for $\Delta e'$, and we use normal coordinates:     
$$
\frac{1}{4} \Delta e' 
=
-\sum_\a \nab \d e'(\zab,\za)
=
=
-\sum_\a \zab.\za.e'~.
$$
Now, 
$\za.\la \duzf,\ov{\duzf}\ra=\la \nab_{\za}\duzf,\ov{\duzf}\ra+\la
\dzuf,\nab_{\za}\ov{\duzf}\ra$. 
The map $f$ is pluriharmonic and therefore
$\nab \df(Z,\bar W)=\nab \df (\bar Z,W) =0$ for all $Z,W$ in
$T^{1,0}\Bm$. Since $\Bm$ is K{\"a}hler,  
$\nab_{\zab}\duzf={\big(\nab_{\zab}\df\big)}^{1,0}$ and hence vanishes
identically on $T^{1,0}\Bm$. It follows that  
$\nab_{\za}\ov{\duzf}=\ov{\nab_{\zab}\duzf}=0$. Thus,
$$
\begin{array}{rcl}
\displaystyle\frac{1}{4} \Delta e' &=& \displaystyle -\sum_\a
\zab.\la \nab_\za\duzf,\ov{\duzf}\ra\\
 &=& \displaystyle -\sum_\a \Big[
\la \nab_{\zab}\nab_\za\duzf,\ov{\duzf}\ra+\la \nab_\za\duzf,\ov{\nab_{\za}\duzf}\ra\Big]\\
 &=& \displaystyle -\frac{1}{2}\|\nab'\duzf\|^2 -\sum_\a 
\la \nab_{\zab}\nab_\za\duzf,\ov{\duzf}\ra~.\\
\end{array}
$$

Therefore, since $(\nab_{\zab}\nab_\za\duzf)(\zb)
=(\nab_{\zab}\nab_\za\duzf)(\zb)-(\nab_{\za}\nab_\zab\duzf)(\zb)$, we have
$$
\begin{array}{rcl}
\ds\sum_\a(\nab_{\zab}\nab_\za\duzf)(\zb) & =  & 
\ds \sum_a R^n(\df(\za),\df(\zab))\duzfb - \sum_a \duzf(R^m(\za,\zab)\zb)\\
& = & \ds\sum_a R^n(\df(\za),\df(\zab))\duzfb+\frac{1}{2}\,\duzf({\rm Ric}^m(\zb))~.
\end{array}
$$
The result follows since the Ricci curvature tensor of $g_m$ is $-\frac{m+1}{2}g_m$. 
\end{demo}  

We also have 

\begin{lemma}\label{dlep}
At each point of $\Bm$ where $e'\not=0$, resp. $e''\not=0$,
$$\frac{1}{4}\Delta\log e'=-\alpha'-2\,\frac{R'}{e'}+\frac{m+1}{4}~,
\mbox{ resp.}\;\;\;
\frac{1}{4}\Delta\log e''=-\alpha''-2\,\frac{R''}{e''}+\frac{m+1}{4}~,$$
where
$$\alpha'=\frac{1}{2e'^2}\Big(\|\nabla'\partial^{1,0}f\|^2 e'-\|\la
\nabla'\partial^{1,0}f,\ov{\partial^{1,0}f}\ra\|^2\Big)$$ 
 and 
$$\alpha''=\frac{1}{2e''^2}\Big(\|\nabla'\partial^{0,1}f\|^2 e''-\|\la
\nabla'\partial^{0,1}f,\ov{\partial^{0,1}f}\ra\|^2\Big)$$ 
are both nonnegative by Cauchy-Schwarz's inequality.
\end{lemma}

\begin{demo}
Again, we make this (easy) computation only for $\frac{1}{4}\Delta\log
e'=\frac{1}{4e'}\Delta e'+\frac{1}{4{e'}^2}\|{\rm d}e'\|^2$. 
Now,
$$
\frac{1}{4}\|{\rm d}e'\|^2 
= 
\sum_\a {\rm d}e'(\za){\rm d}e'(\zab)
= 
\sum_\a \la\nab_\za\duzf,\ov{\duzf}\ra\la\duzf,\nab_\zab\ov{\duzf}\ra
= 
\frac{1}{2} \|\la\nab'\duzf,\ov{\duzf}\ra\|^2
$$
where $\la\nab'\duzf,\ov{\duzf}\ra$ denotes the 1-form on $T^{1,0}\Bm$ given
by 
$
\la\nab'\duzf,\ov{\duzf}\ra(z)=\la\nab_z\duzf,\ov{\duzf}\ra
=2\sum_\b g_n\big((\nab_z\duzf)(\zb),\ov{\duzf}(\zb)\big).
$
\end{demo}

Finally, easy computations show that:

\begin{lemma}\label{rprs}
$$
R'=\frac{1}{2}\sum_{\alpha,\beta}\Big[\big|g_n(\duzf(\za),\ov{\duzf(\zb)})\big|^2
-\big|g_n(\duzf(\zb),\ov{\duzf(\zab)})\big|^2\Big]+\frac{1}{8}e'(e'-e'')~,
$$
$$
R''=\frac{1}{2}\sum_{\alpha,\beta}\Big[\big|g_n(\duzf(\zab),\ov{\duzf(\zbb)})\big|^2
-\big|g_n(\duzf(\za),\ov{\duzf(\zbb)})\big|^2\Big]+\frac{1}{8}e''(e''-e')~.
$$
\end{lemma}

\begin{rema}
  In the sequel, we shall use the fact that all the functions involved in
  those three lemmas are well defined on $M=\Gamma\backslash\Bm$.
\end{rema}

\section{Rigidity of representations of lattices of ${\rm PU}(m,1)$ into ${\rm PU}(n,1)$ }\label{rigid}

\subsection{Burger-Iozzi invariant}\label{invariant}\hfill

We again assume that $m\geq 2$.

Let $\Gamma$ be a torsion-free lattice in ${\rm PU}(m,1)$, and let
$\rho:\Gamma\fd{\rm PU}(n,1)$ be a homomorphism. 
M.~Burger and A.~Iozzi assign to $\rho$ an
invariant which can be defined as follows (see \cite{BI01}). 

Take any $\rho$-equivariant map $f:\Bm\fd\Bn$ and consider the pull-back
$f^\star\omega_n$ of the Kähler form $\omega_n$ of $\Bn$. Note that we can
consider $f^\star\omega_n$ as a 2-form on $M=\Gamma\backslash\Bm$. The 
de Rham cohomology class $[f^\star\omega_n]\in H^2_{DR}(M)$ defined by
$f^\star\omega_n$ is independent of the choice of the equivariant map $f$ since
all such maps are homotopic, and therefore we call it  $[\rho^\star\omega_n]$.

Now, Burger and Iozzi remark that the class $[\rho^\star\omega_n]$ is in the image
of the natural comparison map from the $L^2$-cohomology group $H^2_{(2)}(M)$ of $M$ to
the de Rham cohomology group $H^2_{DR}(M)$. Since $m\geq 2$, the comparison map is
injective (see \cite{Zu82}; the arithmeticity of the lattice
$\Gamma$ is not necessary for the result in the present case), this yields a
well-defined $L^2$-cohomology class, denoted by ${[\rho^\star\omega_n]}_{(2)}$, and
they define (in a slightly different form) 
$$
\tau(\rho):=\frac{1}{2m}\int_M \la \rho^\star\omega_n,\omega_m\ra dV_m~,
$$
where $\rho^\star\omega_n$ is any $L^2$-form representing
${[\rho^\star\omega_n]}_{(2)}$ (observe that, because $\om_m$ is parallel,
$\tau(\rho)$ depends only on ${[\rho^\star\omega_n]}_{(2)}$, hence on $\rho$).  

\begin{rema} In complex dimension 1 and for a {\it uniform} lattice $\G$,
$\tau(\rho)$ is also well-defined and coincides
with the classical Toledo invariant (cf.~\cite{To89}).
\end{rema}


The main result of~\cite{BI01} then reads:

\begin{theo}\label{burgio}
Under the above assumptions,
$$\big|\tau(\rho)\big|\leq {\rm Vol}(M)~.$$ 
Moreover, equality holds if and only if
there exists a $\rho$-equivariant totally geodesic
  isometric embedding $\Bm\fd\Bn$.
\end{theo}

Let $\G$ be a
torsion-free lattice in $\SU(m,1)$. Via the natural inclusion of $\SU(m,1)$  
into $\PU(n,1)$ for $n>m$, we obtain the so-called standard representation of
$\G$ into $\PU(n,1)$ (which is of course $\C$-Fuchsian). Theorem~\ref{burgio}
then implies:  

\begin{coro}\label{deform}
Let $\G$ be a torsion-free lattice in $\SU(m,1)$, $m\geq 2$. Then any
deformation of the standard representation of $\Gamma$ into ${\rm PU}(n,1)$
($n>m$) is also $\C$-Fuchsian. 
\end{coro}

\begin{demo} Since $\G$ is torsion-free, it projects isomorphically into
  $\PU(m,1)$. We can therefore consider the standard representation of $\G$
  as a representation of a lattice of $\PU(m,1)$ and apply Theorem~\ref{burgio}.
Now, when seen as a cohomology class in $H^2_{DR}(M)$,
  $[\rho^\star\omega_n]$ 
is a characteristic class of the principal ${\rm PU}(n,1)$-bundle over $M$
associated to $\rho$ and so, it is constant on connected components of ${\rm
  Hom}(\Gamma,{\rm PU}(n,1))$. On the other hand, $\big|\tau(\rho)\big|={\rm
  Vol}(M)$ holds when $\rho$ is the standard
representation of $\G$ into ${\rm PU}(n,1)$, hence the result.
\end{demo}

The main tool in~\cite{To89} and~\cite{BI01} is bounded cohomology. 
Corlette in~\cite{Co88} was the first to obtain  Corollary~\ref{deform} 
for $m\geq 2$ and $\G$ cocompact. He worked with an invariant similar to
$\tau(\rho)$, the {\it volume} of $\rho$ (see the introduction and the remark at
the end of Section~\ref{proof4}), and he
used harmonic maps techniques to obtain a result equivalent 
to Theorem~\ref{burgio}.   

We will now show 
that Theorem~\ref{burgio} 
is a consequence of our results on harmonic maps.

\subsection{Proof of Theorem~\ref{burgio}}\label{proof4}\hfill

We will first prove that the invariant is well-defined, that is, that the
de Rham cohomology class $[\rho^\star\omega_n]$ can be represented by an
$L^2$-form. We will actually show that there exists a compactly supported form
in the class  $[\rho^\star\omega_n]$. This will be done by choosing for each
cusp $C$ of $M$ a Kähler potential for $\omega_n$ invariant by the image of the
fundamental group of $C$. 

Let $f$ be a $\rho$-equivariant map from $\Bm$ to $\Bn$. Let 
$C$ be a cusp of $M$ and $\gamma_0$ a generator of the center of
$\Gamma_{C}$. Depending on the type of $\rho(\gamma_0)$, the whole image
$\rho(\Gamma_C)$ of $\Gamma_C$ will fix a point, either in $\Bn$ or on the
boundary at infinity of $\Bn$. 

If $\rho(\gamma_0)$ is parabolic, we can assume that its fixed point is
$\infty$ in the Siegel model of complex hyperbolic space. We then 
have $\omega_n=-{\rm d}{\rm d}^c t'$ where 
$t'=\log(2{\rm Re}(w')-\la\la z',z'\ra\ra)$
(see section 1). As can be easily checked, the 1-form $\varsigma:=-{\rm d}^c t'$ is
invariant by the stabilizer in $\PU(n,1)$ of the fixed point.
Hence $f^\star\varsigma$ can be seen as a 1-form in the cusp $C$ and we 
have, on $C$, that $f^\star\omega_n={\rm d}f^\star\varsigma$. 

If $\rho(\gamma_0)$ is elliptic, then $\rho(\Gamma_C)$ fixes a point in
$\Bn$. We can assume that the fixed point is $0$ in the ball model of
complex hyperbolic space, so that $\rho(\Gamma_C)\subset{\rm U}(m)$. In
this case, we can write $\omega_n={\rm d}\varsigma$ with $\varsigma :=- {\rm
  d}^c\log(1-\la\la z,z\ra\ra)$ where in this case $z$ is a point in the
unit ball of ${\mathbb C}^n$ and $\la\la\,,\,\ra\ra$ denotes
the standard Hermitian product on ${\mathbb C}^n$. Again, this 1-form is
invariant by ${\rm U}(m)$ and therefore we have $f^\star\omega_n={\rm
  d}f^\star\varsigma$ on $C$. 

Repeating this operation in each cusp $C_i$ of $M$, we obtain that
there exists a 1-form $\varsigma_i$ on $C_i$ such that
$f^\star\omega_n={\rm d}f^\star\varsigma_i$  on $C_i$. Therefore, if $\chi$
is a function on $M$ identically equal to 0 on the compact part of $M$
and to 1 far enough in the cusps of $M$, the form $f^\star\omega_n-{\rm
  d}(\sum_i\chi f^\star\varsigma_i)$ is cohomologous to $f^\star\omega_n$
and has compact support, hence is $L^2$. 

We are now ready to prove the theorem. We begin with the:

\begin{lemma}
Suppose the representation $\rho$ is not reductive. Then $\tau(\rho)=0$.
\end{lemma} 

\begin{demo}
Let $f$ be
any $\rho$-equivariant map from $\Bm$ to $\Bn$. Since
the image $\rho(\Gamma)$ fixes a point at infinity in $\Bn$, $\rho(\Gamma)$
is a subgroup of the stabilizer of this fixed point in $\PU(n,1)$ and
therefore, the form $f^\star\varsigma$ on $\Bm$ defined above goes down to
a form defined on the whole quotient $M$ such that $f^\star\omega_n={\rm
  d}f^\star\varsigma$. Hence the de Rham cohomology class of
$f^\star\omega_n$ is zero and the invariant $\tau(\rho)$ vanishes.  
\end{demo}

We may now assume that $\rho$ is reductive. Theorems \ref{existharmonic} and
\ref{finiteenergy} then guarantee the existence of a finite energy
$\rho$-equivariant harmonic map $f:\Bm\fd\Bn$. 

\medskip

We first prove that we can use $f$ to compute $\tau(\rho)$. This follows from
the 
\begin{lemma}\label{canusef}
The 2-form $f^\star\omega_n$ is $L^2$.  
\end{lemma}

\begin{demo}
If $f$ is (anti)holomorphic, it follows from the
  generalization of the Schwarz-Pick lemma (see~\cite{Ko70})  
  that in an obvious sense, $f^\star g_n\leq g_m$. An easy computation then
  shows that $\|f^\star\omega_n\|\leq 2m$ everywhere on $M$: $f^\star\omega_n$
  is a $L^2$-form. 

Assume now that $f$ is not (anti)holomorphic. Then Proposition \ref{holom}
implies that ${\rm rk}_\R f\leq 2$.
When ${\rm rk}_\R\, f\leq 1$, $f^\star\om_n=0$, whereas ${\rm rk}_\R\,
{\rm d}_xf=2$ (or ${\rm rk}_\C\, {\rm d}^\C_xf=2$) is equivalent to 
$$
{\rm dim}_\C\,{\rm d}^\C_x f(T^{1,0}\Bm)=1\ {\rm and}\ {\rm d}^\C_x
f(T^{1,0}\Bm)\cap\overline{{\rm d}^\C_x f(T^{1,0}\Bm)}=\{0\}~.
$$ 
But, if ${\rm d}^\C_x f(T^{1,0}\Bm)$ contains no real vectors, it contains no
purely imaginary vectors, thus ${\rm d}_x f(X)=0$ if and only if ${\rm d}^\C_x
f(X-iJX)=0$. This means that ${\rm Ker}\, {\rm d}_x f$ is $J$-invariant, hence,
we may choose an orthonormal basis $(e_i)_{1\leq i\leq 2m}=(e_i,Je_i)_{1\leq
  i\leq m}$ of $T_x\Bm$ in which computations give 
$$
\la f^\star\om_n,f^\star\om_n\ra = \sum_{i,j=1}^{2m} f^\star\om_n(e_i,e_j)^2
= 2 f^\star\om_n(e_1,Je_1)^2 = 2 (e'-e'')^2~.
$$

Thus, if ${\rm rk}_\R f=2$, we have $\|f^\star\om_n\|^2=2(e'-e'')^2$ on
$\Bm$. Using the fact that ${\rm rk}_\R f=2$ in the formulae of
Lemma~\ref{rprs}, we find that $R'\geq\frac{1}{4}e'(e'-e'')$ and 
$R''\geq\frac{1}{4}e''(e''-e')$ 
so that $(e'-e'')^2\leq 4(R'+R'')$. Now, adding the two equalities in
Lemma~\ref{dep}, we get   
$$
2(e'-e'')^2\leq 8(R'+R'')=-\Delta e-2\|\nabla'{\rm d}^\C f\|^2+(m+1)\, e~.
$$
Since $f$ has finite energy, we conclude that $f^\star\om_n$ is $L^2$ (see the
proof of Lemma~\ref{nabladfl2} and the beginning of section~\ref{lemmestechniques}). 
\end{demo}

We now prove the inequality $|\tau(\rho)|\leq{\rm Vol}(M)$.

The above lemma implies that the Burger-Iozzi invariant of $\rho$ is given by 
$$
\tau(\rho)=\frac{1}{2m}\int_M \la f^\star\omega_n,\omega_m\ra
dV_m=\frac{1}{m}\int_M (e'-e'')dV_m~, 
$$  
since $\la f^\star\omega_n,\omega_m\ra=2(e'-e'')$.

Therefore, if $f$ is a complex rank $r$ holomorphic map, the Schwarz-Pick lemma
implies that $0\leq\la f^\star\om_n,\om_m\ra\leq 2r$ at each point of $M$,
whereas if $f$ is a rank $r$ antiholomorphic map, then $-2r\leq\la
f^\star\om_n,\om_m\ra\leq 0$. Integrating these inequalities yields the result.

If $f$ is neither holomophic nor antiholomorphic, then we know from
Proposition~\ref{holom} that ${\rm rk}_\R f\leq 2$. We will prove that in this
case,
$$
\Big|\int_M\la f^\star\omega_n,\omega_m\ra dV_m\Big|\,\leq\,(m+1){\rm Vol}(M)~.
$$

As in the
proof of Lemma~\ref{canusef}, we have  $R'\geq\frac{1}{4}e'(e'-e'')$ and 
$R''\geq\frac{1}{4}e''(e''-e')$. Lemma~\ref{dlep} then implies that 
$$
\Delta\log e'\leq 2(e''-e')+m+1~,
$$
resp.
$$
\Delta\log e''\leq 2(e'-e'')+m+1~,
$$
at each point where $e'\not=0$, resp. $e''\not=0$.

Let $\varepsilon >0$ be a regular value of $e'$. We set $M_\varepsilon=\{x\in
M,\,e'(x)>\varepsilon\}$, and we introduce the cut-off function $\eta_R$
defined in the proof of Proposition~\ref{bochnerformula}.

By Green's formula we have
$$
\int_{M_\varepsilon}\Big[\eta_R\Delta\log e'-\la {\rm grad}\,\eta_R,{\rm
  grad}\, \log e'\ra\Big]dV_m=\int_{\partial M_\varepsilon} \eta_R \la\nu,{\rm
  grad}\, \log e'\ra\,dA\geq 0~.
$$ 
The latter is nonnegative because $\nu$ is the inward unit vector field along
$\partial M_\varepsilon$ which is pointwise orthogonal to $\partial
M_\varepsilon$. From 
$$
|\la {\rm grad}\,\eta_R,{\rm grad}\, \log e'\ra|\leq \sqrt2\|{\rm
  d}\eta_R\|\frac{\|\nabla'\partial^{1,0}f\|}{\sqrt{e'}}\leq \sqrt2\|{\rm
  d}\eta_R\|\frac{\|\nabla'\partial^{1,0}f\|}{\sqrt\varepsilon}\ \ {\rm on}\
M_\varepsilon~,
$$ 
we get
$$
\int_{M_\varepsilon}\eta_R\,\Delta\log e'\,dV_m\geq
-\sqrt{\frac{2}{\varepsilon}}\frac{C}{R} \int_{M}\|\nabla'\partial^{1,0}f\|
dV_m\geq -\sqrt{\frac{2}{\varepsilon}}\frac{C}{R}\Big[{\rm Vol}
(M)\int_{M}\|\nabla'\partial^{1,0}f\|^2 dV_m\Big]^{1/2}~.
$$
Using the fact that $\|\nabla'\partial^{1,0}f\|\in L^2(M)$,
we obtain that, for some constant $C'$ and for all $R>1$, 
$$
\int_{M_\varepsilon}\eta_R\Delta\log
e'\,dV_m\geq -\sqrt{\frac{2}{\varepsilon}}\frac{C'}{R}
$$ 
and so
$$
\int_{M_\varepsilon} \Bigl(2(e''-e')+m+1\Bigr)dV_m=\lim_{R\pfd
  +\infty}\int_{M_\varepsilon} \eta_R\Bigl(2(e''-e')+m+1\Bigr)dV_m\geq 0~.
$$ 
The subset $\{x\in M,\,e'(x)=0\}$ has zero measure since $\duzf$ is a
holomorphic section of ${\rm Hom}(T^{1,0}\Bm,f^\star T^{1,0}\Bn)$. Hence, 
letting $\varepsilon\rightarrow 0$, we conclude that
$$
\int_M \la f^\star\omega_n,\omega_m\ra dV_m=\int_M 2(e'-e'')dV_m\leq (m+1)\,{\rm
  Vol}(M)~.
$$

Integrating $\Delta\log e''$ in the same way, we get the required inequality.

\medskip

Assume now that equality holds in the theorem: 
$|\int_M\la f^\star\omega_n,\omega_m\ra dV_m|=2m{\rm Vol}(M)$.

Since $m\geq 2$, it follows from the proof above 
that $f$ is a complex rank $m$ holomorphic or antiholomorphic map. Since the
inequality $|\la f^\star\om_n,\om_m\ra|\leq 2m$ is therefore true 
pointwise, the global equality implies that  $|\la f^\star\om_n,\om_m\ra|=2m$ 
everywhere on $M$: $f$ is an isometry. On the other hand, the
Bochner-type formula ($\diamondsuit$) with $Q=I$ reads:
$$
\int_M\|\nab\df\|^2dV_m = \frac{1}{2}\int_M\Big({\rm Scal} (f^\star
R^n)-\frac{1}{2m}\|\df\|^2{\rm Scal} (R^m)\Big)dV_m~. 
$$
Since $f$ is an isometry, $\|\df\|^2=2m$ and $f^\star R^n=R^m$.
Therefore $\nab\df=0$, namely, $f$ is totally geodesic and we are done.


\begin{rema}
It should be noted that if one is interested only in proving
Corollary~\ref{deform}, it is actually possible to 
define another invariant, that one could call the
$L^2$-volume of the representation, in the following way: just take the
$m$-th exterior power of any $L^2$-form representing the $L^2$-cohomology
class $[\rho^\star\omega_n]_{(2)}$ and integrate it over $M$. The
so-obtained number is independent of the choice of the $L^2$-representative
and therefore depends only on $\rho$. One can then prove
Theorem~\ref{burgio} with $\tau(\rho)$ replaced by this (suitably
normalized) $L^2$-volume. The proof is in fact easier since this
invariant will be zero if the real rank of the harmonic map is less than
$2m$. Therefore one does not need to deal with non (anti)holomorphic maps.

Nevertheless, there exist representations $\rho$ with zero $L^2$-volume
but $\tau(\rho)\neq 0$ and one can hope to be able to get informations on 
these representations from the Burger-Iozzi invariant that the
volume would not give. An example of such a representation was given by
Livn\'e in his thesis (\cite{Livne}, see also \cite{Kapovich}). He
constructed a (closed) complex hyperbolic manifold $M$ of complex dimension 2
and a surjective holomorphic map $f$ from $M$ to a (closed) Riemann surface
$\Sigma$ such that the induced map on the fundamental groups is surjective.
This gives a representation of the lattice $\pi_1(M)\subset\PU(2,1)$ into
$\PU(1,1)$. Of course the volume of this representation is zero. But its
Burger-Iozzi invariant does not vanish since the pull-back
$f^\star\omega_1$ is a semi-positive $(1,1)$-form on $M$ which is not
identically zero.  
\end{rema}

\section{The case of non-uniform lattices of ${\rm PU}(1,1)$}\label{dim1}

In this section we want to extend the previous results to the 1-dimensional
case, namely the case of non-uniform lattices of ${\rm PU}(1,1)$. 
We remark that if $\rho$ is a representation of such a lattice into
$\PU(n,1)$, the Burger-Iozzi invariant of $\rho$ is not defined
since $H^2_{DR}(M)=0$ and the comparison
map $H^2_{(2)}(M)\fd H^2_{DR}(M)$ is of course not injective
anymore. As we mentionned in the introduction, Corollary~\ref{deform} fails in
this case. Indeed, Gusevskii and Parker prove in~\cite{GP00} that there exist 
lattices in ${\rm PU}(1,1)$ admitting quasi-Fuchsian deformations into
${\rm PU}(2,1)$.   

Geometrically, a torsion-free lattice $\G<\PU(1,1)$ is the fundamental group of
the complete hyperbolic surface of finite volume $M=\G\backslash\B^1$. It turns
out that we don't need a Riemannian structure on $M$ to define an invariant
associated to representations of its fundamental group into $\PU(1,1)$. We will
therefore  work in the more general setting of fundamental 
groups of orientable surfaces of finite topological type. 

Let $M$ be the open surface obtained by removing a finite number of points 
$m_1,\ldots,m_p$, called punctures, from a closed orientable surface
$\ov{M}$ of genus $g$. We will call such an $M$ a $p$-times
punctured closed orientable surface of genus $g$. We assume that the Euler
characteristic 
$\chi(M)=2-2g-p$ of $M$ is negative. Let $\G=\pi_1(M)$ be the fundamental
group of $M$ and $\pi:\wt M\fd M$ be the universal cover of $M$. 

Loops going once around the puncture $m_i$ in the direction prescribed by
the orientation of $\ov M$ correspond to a conjugacy class $c_i$ of
elements of $\G$. The elements of the conjugacy classes $c_i$ are called
{\it peripheral}. For each $i$, choose $\g_i\in c_i$ and denote by
$\la\g_i\ra$ the cyclic subgroup generated by $\g_i$. 
There exist small disjoint open topological discs $D_i\subset \ov M$ around
each $m_i$ and disjoint open simply-connected sets $U_i$ in $\wt M$,
precisely invariant under $\la\g_i\ra$ (meaning that $\g_i U_i=U_i$
and $\g U_i\cap U_i=\emptyset$ if $\g\notin\la\g_i\ra$), such that the
punctured discs $D_i^\star:=D_i\backslash\{m_i\}\subset M$ are given 
by $D_i^\star=\la\g_i\ra\backslash U_i$ (This can for example be seen by
uniformizing $M$ as a finite volume hyperbolic surface, and then choosing
precisely invariant horospherical neighbourhoods of the parabolic fixed
points of $\G$).  

Let $\rho$ be a homomorphism from $\G$ to  $\PU(n,1)$. At the beginning of
section~\ref{proof4} (there, $\G$ was a lattice in $\PU(m,1)$ for $m\geq
2$), we saw how to find a compactly supported 2-form in the de Rham
cohomology class $[\rho^\star\om_n]$. Here, as we said, this class is zero
but we shall in the same way associate to $\rho$ a class in the cohomology
with compact support. 

For each $i$, choose a fixed point $\xi_i$ of $\rho(\g_i)$, in  $\Bn$ if
$\rho(\g_i)$ is elliptic, else in $\partial_\infty\Bn$, and then a Kähler
potential $\psi_i$ of $\om_n$, invariant by the stabilizer in $\PU(n,1)$ of
$\xi_i$ if  $\xi_i\in\Bn$ or by the stabilizer of the horospheres centered
at $\xi_i$ if $\xi_i\in\partial_\infty\Bn$. If $\xi_i\in\Bn$, we can assume that
$\xi_i=0$ in the ball model of $\Bn$ and take $\psi_i=\log(1-\la\la z,z
\ra\ra)$ (here $\la\la\, ,\,\ra\ra$ denotes the standard Hermitian form on
${\mathbb C}^n$). If $\xi_i\in\partial_\infty\Bn$, we can take $\psi_i=\log (2{\rm
  Re}(w')-\la\la z',z'\ra\ra)$, where $(w',z')\in \C\times\C^{n-1}$ are horospherical 
coordinates centered at $\xi_i$. Note that, up to an additive
constant, these potentials are unique. The potential $\psi_i$ is invariant by
$\rho(\g_i)$ only if $\rho(\g_i)$ is elliptic or parabolic but the 1-form
$\varsigma_i:=-\d^c\psi_i$, which satisfies $\d\varsigma_i=\om_n$, is always
invariant by $\rho(\g_i)$. 

Given a $\rho$-equivariant map $f:\wt M\longrightarrow\Bn$, we can
pull-back the Kähler form $\om_n$ and the forms $\varsigma_i$ to get a
2-form $f^\star\om_n$ invariant by $\G$ and 1-forms $f^\star\varsigma_i$
invariant by $\la\g_i\ra$. We can therefore consider $f^\star\om_n$ as a
2-form on $M$ and  each $f^\star\varsigma_i$ as a 1-form on the punctured
disc $D_i^\star$ (by restricting it first to $U_i$). If now $\chi$ is a
function identically equal to 0 outside the $D_i^\star$'s and to 1 in small
neighbourhoods of the punctures, we get a compactly supported 2-form
$f^\star\omega_n-{\rm d}(\sum_i\chi f^\star\varsigma_i)$ on $M$. This
yields a class  $[f^\star\omega_n-{\rm d}(\sum_i\chi
f^\star\varsigma_i)]_c$  in the second cohomology group with compact
support $H^2_c(M)$.  

\begin{propdef}\label{tau}
  The class $[f^\star\omega_n-{\rm d}(\sum_i\chi f^\star\varsigma_i)]_c$
  depends only on the representation $\rho$, and will therefore be denoted
  by $[\rho^\star\omega_n]_c$. Moreover, we set
$$
\tau(\rho)=\int_M[\rho^\star\omega_n]_c~.
$$  
\end{propdef}

\begin{demo}
This class is clearly independent of the cut-off function $\chi$, and
therefore also of the choice of the punctured discs $D_i^\star$ where the
$f^\star\varsigma_i$'s are defined.       

Now, let $f_1$ and $f_2$ be two $\rho$-equivariant maps $\wt
M\longrightarrow\Bn$. The map 
$f_2$ is homotopic to a $\rho$-equivariant map $f_3$ such that, when seen
as sections of the $\Bn$-bundle on $M$ associated to $\rho$, $f_3=f_1$ on the
set $\{\chi<1\}$ and $f_3=f_2$ close enough to the punctures. Then there
exists a compactly supported 1-form $\alpha$ such that
$f_2^\star\omega_n=f_3^\star\omega_n+{\rm d}\alpha$. Hence we have
$$
f_2^\star\omega_n-{\rm d}\Big(\sum_i\chi f_2^\star\varsigma_i\Big)=
f_3^\star\omega_n-{\rm d}\Big(\sum_i\chi f_3^\star\varsigma_i\Big)
+{\rm d}\Big(\sum_i\chi (f_3^\star\varsigma_i- f_2^\star\varsigma_i)\Big)+{\rm d}\alpha~.\\
$$
But 
$$
f_3^\star\omega_n-{\rm d}\Big(\sum_i\chi f_3^\star\varsigma_i\Big)=
\left\{
\begin{array}{l}
f_1^\star\omega_n-{\rm d}(\sum_i\chi f_1^\star\varsigma_i)\mbox{ on
  }\{\chi<1\}\\
0 \mbox{ on  }\{\chi=1\}
\end{array}
\right.
= f_1^\star\omega_n-{\rm d}\Big(\sum_i\chi f_1^\star\varsigma_i\Big)\mbox{ on
  } M~.
$$
Therefore 
$$
f_2^\star\omega_n-{\rm d}\Big(\sum_i\chi f_2^\star\varsigma_i\Big)=
f_1^\star\omega_n-{\rm d}\Big(\sum_i\chi f_1^\star\varsigma_i\Big)
+{\rm d}\Big[\sum_i\chi (f_3^\star\varsigma_i-
f_2^\star\varsigma_i)+\alpha\Big]~.\\ 
$$
From the definition of $f_3$, the 1-form inside the brackets is compactly
supported.   
Hence $[f^\star\omega_n-{\rm d}(\sum_i\chi f^\star\varsigma_i)]_c$
does  not depend on the $\rho$-equivariant map $f$.  

If $\rho(\g_j)$ is elliptic for some $j$, it might fix more than one point
in $\Bn$. We therefore have to check that choosing another fixed point, say
$\xi_j'$, to define the  Kähler potential does not change the class. Let
$\psi_j'$ be the Kähler potential associated to $\xi_j'$ and
$\varsigma_j'=-\d^c\psi_j'$ the corresponding 1-form. To compute our class, we
can choose the $\rho$-equivariant map $f:\wt M\fd\Bn$ to be constant equal to
$\xi_j$ in $U_j$ so that $f^\star\varsigma_j=f^\star\varsigma_j'=0$ on 
$D_j^\star$. Therefore  
$$
f^\star\omega_n-{\rm d}\Big(\sum_i\chi f^\star\varsigma_i\Big)
= f^\star\omega_n-\d\Big(\sum_{i\neq j}\chi f^\star\varsigma_i\Big)~,
$$
and the cohomology class is not affected. 

In the same way, if $\rho(\g_j)$ is hyperbolic for some $j$, we must show
that we can choose any of the two fixed points $\xi_j$ and $\xi_j'$ of 
$\rho(\g_j)$ indifferently. In this case, we can arrange that the
equivariant map $f$ maps $U_j$ to the axis of $\rho(\g_j)$. If we take the
potential $\psi_j$ associated to $\xi_j$ (resp. $\xi_j'$) to define
$\varsigma_j$ then $\varsigma_j=-\d^c t=\d t\circ J$ in horospherical 
coordinates $(z,v,t)$ chosen so that $\xi_j=\infty$ and $\xi_j'=0$
(resp. $\xi_j'=\infty$ and $\xi_j=0$). Since in these coordinates the axis
of $\rho(\g_j)$ is the set $\{z=0,v=0\}$, $f^\star\varsigma_j=0$ on $D_j^\star$.  
Again, 
$$
f^\star\omega_n-{\rm d}\Big(\sum_i\chi f^\star\varsigma_i\Big)
= f^\star\omega_n-\d\Big(\sum_{i\neq j}\chi f^\star\varsigma_i\Big)~.
$$

Finally, it is easy to check that a different choice of the peripheral
elements $\g_i$ (and hence of the $U_i$'s) gives the same cohomology
class. Indeed, let $\g_j$ and $\g_j'=\g\g_j\g^{-1}$ be two elements of the
conjugacy class $c_i$. We denote with primes the objects associated to
$\g_j'$ (for example, $\psi_j'$ is the Kähler potential associated to a
fixed point $\xi_j'$ of $\rho(\g_j')$). If $\rho(\g_j)$, and hence
$\rho(\g_j')$, is elliptic or hyperbolic, we can choose as above the
equivariant map $f$ so that $f^\star\varsigma_i=0$ on $U_i$ and
$f^\star\varsigma_i'=0$ on $U_i'=\g U_i$. Thus we can assume that
$\rho(\g_j)$ and $\rho(\g_j')$ are both parabolic. But then
$\xi_j'=\rho(\g)\xi_j$, $\psi_j'=\psi_j\circ\rho(\g^{-1})$ and
$f^\star\varsigma_j'=(\g^{-1})^\star f^\star\varsigma_j$ on
$\B^1$. Therefore the restrictions of $f^\star\varsigma_j'$ to $U_j'=\g
U_j$ and of $f^\star\varsigma_j$ to $U_j$ induce the same form on
$D_j^\star=\la \g_j\ra\backslash U_j=\la \g_j'\ra\backslash U_j'$. 
\end{demo}

\medskip

Before stating the main theorem of this section, we need the following
definitions. 


\begin{defi}\label{deftame}
A homomorphism $\rho$ from the fundamental group $\G$ of a $p$-times
punctured closed orientable surface $M$ into $\PU(n,1)$ is called tame if
it maps no peripheral element of $\G$ to a hyperbolic isometry of $\Bn$.   

It is called a
uniformization representation if it is an isomorphism onto a torsion-free
discrete subgroup $\rho(\G)<\PU(n,1)$ such that $M$ and
$\rho(\G)\backslash\Bn$ are diffeomorphic (then, necessarily, $n=1$). 
\end{defi}

\begin{theo}\label{surfaces}
  Let $M$ be a $p$-times punctured closed orientable surface of genus
  $g$. Assume that $\chi(M)=2-2g-p<0$. Let $\Gamma$ be the fundamental group of
  $M$ and $\rho:\Gamma\longrightarrow\PU(n,1)$ be a homomorphism. 
Then
  $|\tau(\rho)|\leq -2\pi\chi(M)$. 
Moreover, equality holds if and only if 

  - the image $\rho(\G)$ of $\G$ stabilizes a totally geodesic copy of $\B^1$
  in $\Bn$, and hence $\rho$ can be seen as a homomorphism from $\G$ into
  $\PU(1,1)$, and

  - $\rho:\G\fd\PU(1,1)$ is a uniformization representation.\\
If moreover $\rho$ is tame and equality holds then $\rho(\G)<\PU(1,1)$ is a lattice. 
\end{theo}

\begin{rema}
(i) The number $\tau(\rho)$, up to sign, depends only on the diffeomorphism
type of the surface $M$ in the following sense. If $\phi:M'\fd M$ is a
diffeomorphism and if we consider the homomorphism
$\rho'=\rho\circ\phi_\star$ of the fundamental group $\G'$ of $M'$ into
$\PU(n,1)$, then it is easily checked that $\tau(\rho')=\pm\tau(\rho)$
depending on whether $\phi$ is orientation preserving or reversing.    

(ii) In particular, if $\rho:\G\fd\G'<\PU(1,1)$ is a uniformization
representation, so that $M$ is diffeomorphic to $M':=\G'\backslash\B^1$,
then $\tau(\rho)=\pm\int_{M'}[\om_1-\d(\sum_i\chi\varsigma_i)]$. Let $M_0'$ be
the convex core of $M'$. $M_0'$ is a finite volume complete hyperbolic
surface whose boundary consists of finitely many disjoint
closed simple geodesics $c_k$ (corresponding to the conjugacy classes of
peripheral elements of $\G$ sent by $\rho$ to hyperbolic isometries of
$\B^1$). We can assume that $\chi=0$ on the boundary of $M_0'$. It is easy
to see (cf. the proof of Proposition~\ref{invharmonic} below) 
that the 1-forms $\varsigma_j$ corresponding to punctures of $M_0'$ are
$L^1$ forms and therefore that 
$$
\int_{M_0'}[\om_1-\d(\sum_j\chi\varsigma_j)]=\int_{M_0'}\om_1={\rm Vol}(M_0')=-2\pi\chi(M_0')=-2\pi\chi(M)~.
$$
Moreover,
$$
\int_{M'\backslash
  M_0'}[\om_1-\d(\sum_k\chi\varsigma_k)]=\sum_k\int_{c_k}\d\varsigma_k=0~,
$$
as can be seen from the proof of Proposition~\ref{tau}. Hence
$\tau(\rho)=\pm 2\pi\chi(M)$ as Theorem~\ref{surfaces} says.  

(iii) Our definition of the invariant $\tau$ makes no reference to
a Riemannian structure on the surface $M$. This means that we
can equip $M$ with any Riemannian metric we want and that our results will
be independent of this particular choice. Now, the uniformization theorem
implies that there exist complete hyperbolic metrics of finite
volume on $M$. If we choose such a metric $g_1$ on $M$, we obtain a
uniformization representation $u:\G\fd\PU(1,1)$ and $M$ is diffeomorphic to
$u(\G)\backslash\B^1$. The composition $\rho\circ u^{-1}$ is a
homomorphism from the torsion-free lattice $u(\G)$ of $\PU(1,1)$ to
$\PU(n,1)$. These identifications are completely transparent for the
invariant $\tau$ and we can therefore assume that $\G$ is a torsion-free
lattice in $\PU(1,1)$ and that $M=\G\backslash\B^1$ with its complete
hyperbolic metric $g_1$ of finite volume. The punctured neighbourhoods
$D^\star_i$ of the punctures $m_i$ will then be seen as cusps of $M$ and
will often be denoted by $C_i$.   

The metric on $M$ allows to talk about $L^2$-cohomology groups and if
we call $[\rho^\star\omega_n]_{(2)}$ the image in $H^2_{(2)}(M)$ of
$[\rho^\star\omega_n]_c$ under the comparison map $H^2_c(M)\longrightarrow
H^2_{(2)}(M)$, we have 
$$
\tau(\rho)=\frac{1}{2}\int_M\la [\rho^\star\omega_n]_{(2)},\omega_1\ra dV_1~,
$$ 
where 
$\om_1$ is the Kähler form of $g_1$. 

(iv) If $\rho$ is not reductive, that is, if $\rho(\G)$ fixes a
  point in $\partial_\infty\Bn$, then $\tau(\rho)=0$. Indeed, 
we can use this fixed point to define as before a 1-form $\varsigma$
invariant by $\rho(\G)$. Then $f^\star\omega_n={\rm d}f^\star\varsigma$ on
$M$. Hence 
$f^\star\omega_n-{\rm d}(\chi f^\star\varsigma) ={\rm d}((1-\chi)
f^\star\varsigma)$, that is $[f^\star\omega_n-{\rm d}(\chi
f^\star\varsigma)]_c=0$ and $\tau(\rho)=0$.

Hence we assume from now on
that all considered representations are reductive.   
\end{rema}

The proof of Theorem~\ref{surfaces} is easier if one deals only with {\it tame}
representations. Since almost all the ideas
are needed in this case, we present it separately in Section~\ref{wtpr}
and we explain in Section~\ref{gc} how to adapt the arguments to handle the
general case.

\subsection{Tame representations}\label{wtpr}\hfill

Note that a representation $\rho$ of a torsion-free lattice $\G$ of $\PU(1,1)$ into
$\PU(n,1)$ is tame if and only if $\rho$ maps no parabolic elements of $\G$ to
hyperbolic elements of ${\rm PU}(n,1)$. The name ``tame'' is 
motivated by the following proposition:

\begin{prop}\label{finenergy}
Let $\Gamma<{\rm PU}(1,1)$ be a torsion-free lattice, and let
$\rho:\Gamma\fd{\rm PU}(n,1)$ be a homomorphism. There exists a
$\rho$-equivariant map $f:\B^1\fd\Bn$ of finite energy if and only if
$\rho$ is tame.
\end{prop}

\begin{demo}
When $\rho$ is tame, we may
easily construct a $\rho$-equivariant map with finite energy. In fact, we
use the same method as in Theorem~\ref{finiteenergy}; since the fundamental
group of each cusp of $M$ is generated by a single parabolic element, the
construction is much simpler. 
Namely, it is sufficient to define the $\rho$-equivariant map on each cusp of
$M$. Let $C$ be a cusp of $M$ and let $\gamma$ be a parabolic element of $\G$
generating $\pi_1(C)$ (via the usual identification). 
We can choose (horospherical) coordinates $(v,t)$ on $\B^1$ such
that $\g(v,t)=(v+a,t)$ and such that $C$ is isometric to the quotient by
$\la \g\ra$ of the subset $D:=[0,a]\times [0,+\infty)$ of $\B^1$ endowed
with the metric $g_1=e^{-2t}dv^2+dt^2$.

If $f:\B^1\fd\Bn$ is any $\rho$-equivariant map, the energy of
$f$ in the cusp is given by 
$$E_C(f)=\frac{1}{2}\int_0^{+\infty}\int_0^a\|{\rm d}f\|^2_{(v,t)}e^{-t}dv\,dt~.$$

When $\rho(\gamma)$ is elliptic, we map $D$ to a fixed point of
$\rho(\gamma)$ and the energy of $f$ in the cusp is zero. 
If $\rho(\gamma)$ is parabolic, let $(z',v',t')$ be adapted horospherical 
coordinates on $\Bn$. We define 
$$
\begin{array}{rcl}
f:[0,a]\times [0,+\infty) & \longrightarrow &\Bn \\
(v,t) & \longmapsto & (\varphi(v),2t)\\
\end{array}
$$
where $\varphi$ is a map from $[0,a]$ into a horosphere $HS'\subset\Bn$
such that $\varphi(a)=\rho(\gamma)\varphi(0)$. 
Computing as in the proof of Theorem~\ref{finiteenergy} ({\it case 1:
  $\rho(\g_0)$ is parabolic}), we get the following estimation of the energy of
$f$ in the cusp: 
$$
E_C(f)\leq \frac{1}{2}\int_0^{+\infty}\int_0^a (4+\|\d\varphi\|^2) e^{-t}d v\,d
t\leq 2a+E(\varphi)<+\infty~.
$$

Suppose now that $\rho(\gamma)$ is hyperbolic and let $f:\B^1\fd\Bn$ be any
$\rho$-equivariant map.  
For every $t$, we denote by $c_t:[0,a]\fd\Bn$ the curve $v\longmapsto f(v,t)$.
Since $\rho(\gamma)$ is hyperbolic, there exists $\delta>0$ such that the
distance between $f(0,t)$ and $f(a,t)=\rho(\gamma)f(0,t)$ is at least
$\delta$, hence the length $l(c_t)$ of $c_t$ is at least $\delta$.
This implies that 
$$
\frac{1}{2}\int_0^a\|{\rm d}f\|^2_{(v,t)}dv\geq e^{2t}
E(c_t)=\frac{e^{2t}}{2}\int_0^a\|{\rm d}c_t\|^2 dv\geq
\frac{e^{2t}}{2a}\bigl(l(c_t))^2\geq \frac{e^{2t}}{2a}\delta^2
$$ 
and
$$
E_C(f)\geq\frac{\delta^2}{2a}\int_0^{+\infty}e^t dt=+\infty~.
$$ 
\end{demo}

The next proposition shows that for a tame representation
$\rho$, the invariant $\tau(\rho)$ can be computed with any finite
energy $\rho$-equivariant map $\B^1\longrightarrow\Bn$.

\begin{prop}\label{invharmonic}
Let $\G<{\rm PU}(1,1)$ be a torsion-free lattice and let $\rho:\G\fd{\rm
  PU}(n,1)$ be a tame homomorphism. Then, 
$$
\tau(\rho)=\frac{1}{2}\int_M \la f^\star\omega_n,\omega_1\ra \,dV_1
$$
for any finite energy $\rho$-equivariant map $f:\B^1\fd\Bn$. 
\end{prop}

\begin{rema}
One could therefore take this as a definition of the invariant for tame
representations. This gives a formulation very similar to the 
classical one for the Toledo invariant of a uniform lattice. 

Note that the energy
finiteness assumption is necessary as the following simple example shows.
Let $\G$ be a torsion-free lattice in ${\rm SU}(1,1)$ generated by two hyperbolic
elements $\gamma_1$ and $\gamma_2$ such that $\gamma_0=[\gamma_1,\gamma_2]$
is parabolic. Then $M=\G\backslash\B^1$ is diffeomorphic to a
once-punctured torus. Let $\rho$ be the inclusion $\G\fd\PU(1,1)$. Keeping
the same notations as in the proof of Proposition~\ref{finenergy}, we can
define $C^\infty$ maps 
$f_\mu:M\fd M$ which equal identity outside the cusp, and given by 
$
f_\mu:[0,a]\times [0,+\infty) 
\longrightarrow 
[0,a]\times \R,~ 
(v,t) 
\longmapsto 
(v,\mu(t))
$
in $D$, where $\mu\in C^\infty([0,+\infty),\R)$. The energy of $f_\mu$ in the cusp is
$
\frac{a}{2}\int_{0}^{+\infty}(e^{t-2\mu(t)}+\mu'(t)^2e^{-t})\,dt$.
Moreover, a formal computation gives:
$\frac{1}{2}\int_C \la f_\mu^\star\omega_1,\omega_1\ra\,dV_1=\int_C
f_\mu^\star\om_1= 
a\int_{0}^{+\infty}\mu'(t)e^{-\mu(t)}\,dt=a(1-\lim_{t\pfd+\infty}e^{-\mu(t)})
$. 
Note that this limit needs not exist and that
for every $c\in[-\infty,\tau(\rho)]$, we may choose
$\mu$ such that $\int_M f_\mu^\star\om_1=c$.
But if the energy of $f_\mu$ is finite, then $\mu(t)\pfd+\infty$ as
$t\pfd+\infty$, and then, $\int_M f_\mu^\star\om_1=\int_M \om_1=\tau(\rho)$. 
\end{rema}

\begin{demo}
Let $f:\B^1\fd \Bn$ be a finite energy $\rho$-equivariant map. Since
$f$ has finite energy, $\la f^\star\omega_n,\omega_1\ra=2(e'(f)-e''(f))$ is
integrable on $M$ and we can write
$$
\tau(\rho)=\frac{1}{2}\int_M \la f^\star\omega_n,\omega_1\ra\, dV_1
-\frac{1}{2}\int_M\la {\rm d}(\sum_i\chi f^\star\varsigma_i),\omega_1\ra\, dV_1~.
$$
The function $\la {\rm d}(\sum_i\chi f^\star\varsigma_i),\omega_1\ra$ is
integrable on $M$ (because so is  $\la f^\star\omega_n,\omega_1\ra$ in the
cusps of $M$). Therefore, if we prove that the 1-form $\sum_i\chi
f^\star\varsigma_i$ is an $L^1$-form, the Stokes formula of~\cite{Gaffney} will
apply and $\int_M\la {\rm d}(\sum_i\chi f^\star\varsigma_i),\omega_1\ra\, dV_1$
will vanish as wanted. 

Now, if the generator $\gamma_i$ of the fundamental group of the cusp $C_i$
is sent by $\rho$ to a parabolic element of $\PU(n,1)$, we know that in
horocyclic coordinates relative to the fixed point of $\rho(\gamma_i)$ we
have $g_n=dt'^2+\varsigma_i^2+4e^{-t'}\la\la dz', dz'\ra\ra$. Therefore,
$$
\|{\rm d}f\|^2={\rm tr}_{g_1} (f^\star g_n)\geq \|f^\star\varsigma_i\|^2~.
$$

If $\rho(\gamma_i)$ is elliptic, we can assume that one of its fixed points is 0 in
the ball model of $\Bn$. In this case, we have
$\varsigma_i=-{\rm d}^c\psi_i$ with $\psi_i=\log(1-\la\la z,
z\ra\ra)$. The metric $g_n$ is then given by 
$$
g_n = 4\frac{\la\la dz, dz\ra\ra}{1-\la\la z,z\ra\ra}+ ({\rm d}\psi_i)^2 
+({\rm d}^c\psi_i)^2~.
$$  
Again, we have $\|f^\star\varsigma_i\|\leq\|{\rm d}f\|$ on $C_i$.

The form $f^\star\varsigma_i$ is therefore $L^1$ on $C_i$ for each $i$. The
proposition follows.  
\end{demo}

\medskip

We are now in position to prove Theorem~\ref{surfaces} in the case of reductive
and tame representations. 

\begin{demo}[of Theorem~\ref{surfaces}]
Since $\rho$ is assumed to be
tame and reductive, Proposition~\ref{finenergy} and
Theorem~\ref{existharmonic} imply the existence of a $\rho$-equivariant harmonic
map  $f:\B^1\fd\Bn$ of finite energy. 

Proceeding exactly as in section~\ref{proof4} yields the inequality
$|\tau(\rho)|\leq{\rm Vol}(M)$.  

Before treating the equality case, we remark that there are particular
representations for which it is immediate that equality can not hold:

\begin{claim}\label{noid}
Suppose that $\rho$ maps the conjugacy classes $c_{q+1}$,\dots, $c_p$ of
peripheral elements 
of $\Gamma$ to ${\rm id}\in\PU(n,1)$ and denote by
$M'$ the surface obtained by removing only the points $m_1$,\dots, $m_q$ from
the closed orientable surface $\overline M$ of genus $g$. Then,
$|\tau(\rho)|\leq\max\bigl(0,-2\pi\chi(M')\bigr)$.
\end{claim}

\begin{demo}
The representation $\rho$ admits a factorisation $\rho=\rho'\circ i_\star$ where
$i_\star:\G\fd\G':=\pi_1(M')$ is the representation induced by the inclusion
$i:M\fd M'$ and $\rho':\G'\fd\PU(n,1)$ is a homomorphism. We may equip $M'$ with
a complex structure such that the inclusion $i$ is holomorphic.

If $\chi(M')<0$, we endow $M'$ with a complete hyperbolic metric of
finite volume. Then,  $|\tau(\rho')|\leq -2\pi\chi(M')$ by the above
argument.
If $f$ is any section of the bundle $M'\times_{\rho'}\B^n$,  
the restriction of $f$ to $M$ can be seen as a section of
the bundle $M\times_{\rho}\B^n$ whose energy is finite on neighbourhoods of
the punctures  $m_{q+1}$,\dots, $m_p$ (since energy
finiteness only depends on the conformal structure of $M$). Then, by
Proposition~\ref{invharmonic}, $\tau(\rho)=\int_M
[f^\star\om_n-\d(\sum_{i=1}^{q}\chi
f^\star\varsigma_i)]=\int_{M'}[f^\star\om_n-\d(\sum_{i=1}^{q}\chi
f^\star\varsigma_i)]=\tau(\rho')$. 

If $\chi(M')=0$, then $M'$ is holomorphically equivalent to either $\C^\star$ or
an elliptic curve. In the first case, $\G'\simeq \Z$ and so $\rho(\G)$ is
generated by a single element. If this element is parabolic or hyperbolic, then
$\rho$ is not reductive and if it is elliptic, there exists a constant
$\rho$-equivariant map. Thus, $\tau(\rho)=0$. In the second case, $\G'$ is
abelian, generated by two elements $\g_1$ and $\g_2$. If $\rho'(\g_1)$ is
parabolic or hyperbolic, by commutation, $\rho'(\g_2)$ must fix the fixed
point(s) of $\rho'(\g_1)$ and $\rho$ is not reductive. If $\rho'(\g_1)$ and
$\rho'(\g_2)$ are elliptic, they must have a common fixed point in $\Bn$, as
they commute. In either case, $\tau(\rho)=0$.

Finally, if $\chi(M')>0$, $M'$ is simply connected and $\tau(\rho)$ is of course zero.
\end{demo}

\medskip

We suppose now that the equality $|\tau(\rho)|=-2\pi\chi(M)$ holds. 

The harmonic map
 $f:\B^1\fd\Bn$ needs not be (anti)holomorphic as in the higher dimensional case, but
 we know that its real rank is 2  
at some point, hence on a dense open subset ${\cal U}$ of $\B^1$ by a result of
\cite{Sa78}, and that in 
the proof of Theorem~\ref{burgio}, one of the inequalities concerning
$R'$ and $R''$ is necessarily an equality on $M$.

Suppose for example that $R'=\frac{1}{4} e'(e'-e'')$ (the case where the
inequality becoming an equality is $R''\geq\frac{1}{4}e''(e''-e')$ is
handled similarly), that is 
$$
\Bigl|g_n\Bigl(\partial_x^{1,0}f(z_1),\ov{\partial_x^{1,0}f(\bar z_1)}\Bigr)\Bigr|^2
=g_n\Bigl(\partial_x^{1,0}f(z_1),\ov{\partial_x^{1,0}f(z_1)}\Bigr)\,
g_n\Bigl(\partial_x^{1,0}f(\bar z_1),\ov{\partial_x^{1,0}f(\bar z_1)}\Bigr)
$$
for every $x\in\B^1$.
 
This has the following simple but very important consequence: for all
$x\in\B^1$, ${\rm d}_x f (T_x\B^1)$
is contained in a complex one-dimensional subspace of $T_{f(x)}\Bn$. To
prove this fact, we
only need to consider points where ${\rm rk}_\R\, {\rm d}_x f=2$. We
remark that a vector subspace $V\subset T_{f(x)}\Bn$ is $J$-invariant if
and only if
$$V^\C=V^{1,0}\oplus V^{0,1}=V^{1,0}\oplus\overline{V^{1,0}}$$
where $V^\C\subset T^\C_{f(x)}\Bn $ is the complexification of $V$ and
$V^{1,0}$ resp. $V^{0,1}$ is the projection of $V^\C$ on
$T_{f(x)}^{1,0}\Bn$ resp. $T_{f(x)}^{0,1}\Bn$. So, if $V={\rm d}_x
f(T_x\B^1)$ is a 2-dimensional real subspace of $T_{f(x)}\Bn$, $V$ is
$J$-invariant if and only if ${\rm dim}_\C\,V^{1,0}=1$. Since $V^{1,0}$
is spanned by $\partial_x^{1,0}f(z_1)$ and $\partial_x^{1,0}f(\bar
z_1)$, the above equality holds if and only if ${\rm d}_x f (T_x\B^1)$
is contained in a complex one-dimensional subspace of $T_{f(x)}\Bn$,
according to the equality part of Cauchy-Schwarz's inequality.

Now, we will use the same trick as A.~G.~Reznikov in his paper
\cite{Re93}. Namely, we  equip the product $\B^1\times\Bn$ with the metric
$h_\varepsilon:=\varepsilon g_1+g_n$ (for a given $\e>0$) and we 
consider the map $\phi:\B^1\fd \B^1\times\Bn$ given by
$\phi(x)=(x,f(x))$. $\phi$ is then a harmonic embedding of $\B^1$ into
$\B^1\times\Bn$. Let $\kappa_{\phi^\star h_\varepsilon}$ be the Gaussian
curvature of the induced metric $\phi^\star h_\varepsilon$ on $\B^1$. Since
$\phi$ is harmonic,  Lemma C.4 of~\cite{Re93} implies that
for all $x\in\B^1$, 
$\kappa_{\phi^\star h_\varepsilon}(x)$  is less than or equal to the value
of the sectional curvature of $h_\varepsilon$ on the 2-plane
$T_{\phi(x)}\phi(\B^1)\subset T_{\phi(x)}(\B^1\times\Bn)$ (which we will call
$\kappa_{h_\varepsilon}(T_{\phi(x)}\phi(\B^1))$).

Since $\phi$ is equivariant
w.r.t. the isometric diagonal action of $\Gamma$  on $\B^1\times\Bn$, these
quantities are well defined over the Riemann surface
$M=\Gamma\backslash\B^1$. Moreover, we have

\begin{claim}\label{gaussbonnet}
For each $\e>0$,
$${\rm Area}_{g_1}(M)=-2\pi\chi(M)=-\int_M  \kappa_{\phi^\star h_\varepsilon}\,
  d V_{\phi^\star h_\varepsilon}$$
where $d V_{\phi^\star h_\varepsilon}$ denotes the volume element of the
pull-back metric $\phi^\star h_\varepsilon$ on $M$.
\end{claim}

\begin{demo}
We shall apply the main theorem of~\cite{Li97} to $\phi^\star
h_\varepsilon$. This metric is clearly complete (because $g_1$ is), and has everywhere non-positive
curvature. We have to prove that

$$\lim_{r\pfd +\infty} \frac{{\rm Area}_{\phi^\star h_\varepsilon}\bigl(B_{\phi^\star h_\varepsilon}(r)\bigr)}{r^2}=0$$
where $B_{\phi^\star h_\varepsilon}(r)$ is the geodesic ball
of radius $r$ (at a fixed point that belongs, for example, to the compact part
of $M$) w.r.t. the metric $\phi^\star h_\varepsilon$. Now, simple computations show that
$$dV_{\phi^\star h_\varepsilon}=\sqrt{\e^2+2\e\, e(f)+\frac{1}{2}\|f^\star\om_n\|^2}\,dV_1~.$$
Since $f$ has finite energy and $f^\star\om_n$ is $L^2$, the metric $\phi^\star
h_\varepsilon$ has in fact finite volume and we obtain the expected limit.
Then, it follows from~\cite{Li97} that $\phi^\star h_\varepsilon$ has finite
total curvature and this implies the result (as it is explained in the same paper).
\end{demo}

On ${\cal U}$, the pull-back $f^\star g_n$ is a Riemannian metric.
Let  $\kappa_{f^\star g_n}$ denote its Gaussian curvature.  Again, since $f$ is
harmonic, for all $x\in{\cal U}$, $\kappa_{f^\star g_n}(x)$ is less than or
equal to  $\kappa_{g_n}(T_{f(x)}f(\B^1))$, the sectional curvature of the
2-plane $T_{f(x)}f(\B^1)\subset T_{f(x)}\Bn$ w.r.t. the metric $g_n$. Since
$T_{f(x)}f(\B^1)$ is a complex line in $T_{f(x)}\Bn$, its sectional
curvature is $-1$. Hence, for all $x\in{\cal U}$, $\kappa_{f^\star
  g_n}(x)\leq -1$.

Assume that at a certain point $p$ of ${\cal U}$, $\kappa_{f^\star
  g_n}(p)$ is stricly less than $-1$. We are going to prove that this can not
happen. For, if it does, there exist $\a>0$ and a small
disc $D$ around $p$ whose closure $\ov D$ is contained in ${\cal U}$
such that for all $x\in D$, $\kappa_{f^\star g_n}(x)\leq -1-\a$. We can
assume that the image of $D$ in $M$ is diffeomorphic to $D$, and we still call
it $D$. Now, if $K$ is any compact subset of ${\cal V}
:=\Gamma\backslash{\cal U}$ containing $D$, we have  
$$
\begin{array}{rcl}
{\ds {\rm Area}_{g_1}(M)}
& \geq &
{\ds -\,\int_K
\kappa_{\phi^\star h_\varepsilon}(x)\,d V_{\phi^\star
  h_\varepsilon}}\\
& \geq & {\ds
-\,\int_D \kappa_{\phi^\star h_\varepsilon}(x)\,d V_{\phi^\star
  h_\varepsilon}
-\,\int_{K\backslash D} \kappa_{h_\varepsilon}(T_{\phi(x)}\phi(\B^1))\,d V_{\phi^\star
  h_\varepsilon}}~,
\end{array}
$$
since on $K$, $\kappa_{\phi^\star
  h_\varepsilon}\leq\kappa_{h_\varepsilon}\leq 0$.

When $\e$ goes to 0, the induced metric $\phi^\star h_\e$ clearly goes to
$f^\star g_n$ on $D$ and therefore, on $D$, the curvature and the volume form
of  $\phi^\star h_\e$ respectively go to the curvature and volume form of
$f^\star g_n$. In the same manner, on the compact set $K$,
$\kappa_{h_\varepsilon}(T_{\phi(x)}\phi(\B^1))$ goes to
$\kappa_{g_n}(T_{f(x)}f(\B^1))$. Again, all involved quantities are well
defined on ${\cal V}$. Hence
$$
{\rm Area}_{g_1}(M)
\,\geq\,
-\,\int_D \kappa_{f^\star g_n}(x)\,d V_{f^\star g_n}
-\,\int_{K\backslash D} \kappa_{g_n}(T_{f(x)}f(\B^1))\,d V_{f^\star
  g_n}~.
$$  
Using $\kappa_{g_n}(T_{f(x)}f(\B^1)) = -1$ and remembering our assumption
on $D$, we obtain
$$
{\rm Area}_{g_1}(M)
\,\geq\, (1+\a)\,{\rm Area}_{f^\star g_n}(D)+{\rm Area}_{f^\star
  g_n}(K\backslash D)\,=\, \a\,{\rm Area}_{f^\star g_n}(D)+{\rm
  Area}_{f^\star g_n}(K)~.
$$
This is true for any compact subset $K$ of ${\cal V}$, hence
$$
{\rm Area}_{g_1}(M)\,\geq\,
\a\,{\rm Area}_{f^\star g_n}(D)+{\rm Area}_{f^\star g_n}({\cal V})~.
$$
Since $f({\cal U})$ is (locally) a 1-dimensional complex submanifold of
$\Bn$, the volume form of
the metric induced on it by $g_n$ is simply the restriction of
$\omega_n$. Therefore the volume form of the metric $f^\star g_n$ on ${\cal
V}$ is $|f^\star\omega_n|=\frac{1}{2}\,|e'-e''|\,\omega_1$. Hence
$$
\begin{array}{rcl}
{\rm Area}_{g_1}(M) & \geq &
\a\,{\rm Area}_{f^\star g_n}(D)+\frac{1}{2}\,\int_{{\cal V}}|e'-e''|\,\omega_1\\
 & = &
\a\,{\rm Area}_{f^\star g_n}(D)+\frac{1}{2}\,\int_{M}|e'-e''|\,\omega_1\\
 & \geq &
\a\,{\rm Area}_{f^\star
  g_n}(D)+\frac{1}{2}\,\Big|\int_{M}(e'-e'')\,\omega_1\Big|~.\\
\end{array}
$$
But we assumed that $\int_{M}(e'-e'')\,\omega_1=2\,{\rm Area}_{g_1}(M)$,
so this gives $e'-e''\geq 0$ on $M$ and the desired contradiction: we
conclude that the curvature of the metric induced by $g_n$ on
$f({\cal U})$ is everywhere $-1$.

If we now take a small open set $U$ in ${\cal U}$ on which $f$ is an
embedding, $f(U)$ is a complex, hence minimal, submanifold  of $\Bn$ whose
sectional curvature equals the restriction of the sectional curvature of
the ambient manifold $\Bn$: $f(U)$ must be totally geodesic and therefore
contained in a complex geodesic. By a result of Sampson (\cite{Sa78}),
this implies that $f$ maps $\B^1$ entirely into this complex geodesic. In particular, $\rho(\Gamma)$
 stabilizes this totally geodesic copy of $\B^1$ and $\rho$ thus induces a tame representation of
$\Gamma$ into $\PU(1,1)$ that we shall still denote by $\rho$ in the sequel. In
the same way, $f$ will be seen as a $\rho$-equivariant map from $\B^1$ into $\B^1$.

\begin{claim}\label{localdiffeo}
The map $f:\B^1\fd\B^1$ is a local diffeomorphism.
\end{claim}

\begin{demo}
We apply the method of Toledo in~\cite{To79} Theorem~4.2, to prove that
$e'-e''>0$ everywhere on $M$. We first go back to the proof of 
Theorem~\ref{burgio}, keeping the same notations. We know that $\Delta\log
e'=2(e''-e')+2$, wherever
$e'>0$. The zeros of $e'$ (if any) are known to be isolated and of finite
order. Suppose that there
exists $x\in M$ such that $e'(x)=0$, and let $r_0,\ \e_0>0$ such that
$e'\geq\e_0$ on  $\partial
{\rm D}(x,r_0)$ (where ${\rm
 D}(x,r_0)$ is the disc of radius $r_0$ about $x$). For any $0<r<r_0$ and
$0<\e<\e_0$, we set $M_{r,\e}=\bigl(M_\e\backslash {\rm D}(x,r_0)\bigr)\cup
\bigl({\rm D}(x,r_0)\backslash{\rm D}(x,r)\bigr)$ and $M'_\e=M_\e\cup{\rm D}(x,r_0)$.

By Green's formula, we have
$$\int_{M_{r,\varepsilon}}\Big[\eta_R\Delta\log e'-\la {\rm grad}\,\eta_R,{\rm
  grad}\, \log e'\ra\Big]dV_1=\int_{\partial M_{r,\varepsilon}} \eta_R \la\nu,{\rm
  grad}\, \log e'\ra\,ds~,$$
and letting $r\pfd 0$, then $R\pfd+\infty$, we get
$$\int_{M'_\varepsilon} \Bigl(2(e''-e')+2\Bigr)dV_1\geq 4\pi q~,$$
where $2q>0$ is the order of the zero of $e'$ at $x$ (see~\cite{To79} and the
proof of Theorem~\ref{burgio}). But, letting $\e\pfd 0$, we get a contradiction
with $\tau(\rho)=-2\pi\chi(M)$.
Finally, we thus have $e'>0$ and $e'-e''\geq 0$ everywhere on $M$.

Now, if $e'(x)=e''(x)$ at some point $x\in M$, there exists a positive
constant $A$ such that 
$$-\Delta\log\frac{e'}{e''}\leq 4 (e'-e'')\leq A\log\frac{e'}{e''}$$
in a neighbourhood of $x$. This implies (see~\cite{SY78} and references
therein) that $\log (e'/e'')$ 
is identically zero on a non empty open subset of $M$, which is impossible, as
the rank of $\d f$ is generically 2.
\end{demo}

We may therefore pull-back the complex structure of the target to obtain a complex
structure on $\B^1$ with respect to which $f$ is holomorphic. In the following
discussion, $M$ will be endowed with the new induced complex structure.

As a consequence of the uniformization theorem, there exists a properly
discontinuous subgroup
$\Gamma'<\PU(1,1)$ such that $M$ is
holomorphically equivalent to $\Gamma'\backslash\B^1$.

\begin{claim}\label{lattice}
The discrete subgroup $\Gamma'<\PU(1,1)$ is a lattice.
\end{claim}

\begin{demo}
Each puncture of $M$ has a neighbourhood
which is holomorphically equivalent to either a punctured disc or an annulus. We
have to show that the latter cannot happen.

Let $m$ be a puncture of $M$ and suppose that $m$ has a neighbourhood $U$
which is holomorphically equivalent to an annulus
$$A_b=\{z\in\C\ ,\ e^{-\pi/b}<|z|<1\}\ ,\ b>0~.$$
In the former complex
structure, we have horocyclic coordinates $(v,t)\in [0,a]\times[0,+\infty)$ on
a fundamental domain of the cusp $C$ corresponding to $m$. Let us denote as before
by $c_t:[0,a]\fd\B^1$ the curve $v\longmapsto f(v,t)$.
Since the harmonic $\rho$-equivariant map $f$ has finite energy, 
there exists a sequence $(t_k)_{k\in\N}$ going to infinity such that
$c_k:=c_{t_k}$ verifies
$$\lim_{k\pfd+\infty}\bigl(l(c_k))^2\leq \lim_{k\pfd+\infty} a
\int_0^a\|{\rm d}c_k\|^2 dv\leq \lim_{k\pfd+\infty} a\,e^{-t_k}\int_0^a\|{\rm
  d}f\|^2_{(v,t_k)}e^{-t_k}dv =0~,
$$
where $l(c_k)$ is the length of $c_k$.

If $\gamma$ is a
peripheral element of $\Gamma$ generating $\pi_1(C)$, $\rho(\gamma)$ is either parabolic
or elliptic, because $\rho$ is tame. In the first case, $f$
induces a holomorphic function $f_U$ from $A_b$ into
$\la\rho(\gamma)\ra\backslash\B^1$ (which is holomorphically equivalent to
a punctured disc). The loops $\alpha_k:[0,a]\fd C$, $v\longmapsto 
(v,t_k)$ define a sequence $(\beta_k)$ of simple loops in $A_b$ (with index
1 with respect to the origin), 
and whose supports converge to the circle $C_b=\{z\in\C\ ,\ |z|=e^{-\pi/b}\}$. 
As $\lim_{k\pfd +\infty} l(c_k)=0$, for every $\e>0$, there exists $k_\e\in\N$
such that $k\geq k_\e$ implies $|f_U\circ\beta_k(t)|\leq \e$ for any $t\in[0,a]$.
 By the maximum principle, $|f_U(z)|\leq\e$ for any $z$ in the connected
 component of $A_b$ lying between the support of $\beta_{k_\e}$ and $C_b$. But the
 holomorphic function $f_U$ admits a Laurent series expansion
 $f_U(z)=\sum_{n\in\Z}a_nz^n$ on $A_b$ and since, for every $e^{-\pi/b}<r<1$,
$$|a_n|=\left|\frac{1}{2i\pi}\int_{|z|=r}\frac
  {f(\zeta)}{\zeta^{n+1}}\,d\zeta\right|\leq e^{n\pi/b}\sup_{|z|=r} |f(z)|~,$$
$f_U$ should be identically zero, which is a contradiction.

If $\rho(\gamma)$ is elliptic, we know from Claim~\ref{noid} that
$\rho(\gamma)\neq{\rm id}$ and thus we may assume that it is a rotation with center
$0\in\B^1$ and angle $\theta\in (0,2\pi)$. 
 We now use the cover map $S\fd A_b$,
$w\longmapsto e^{iz/b}$ where $S$ is the strip $\{w\in\C\ ,\ 0<{\rm Im}\,w<\pi\}$. The
$\rho$-equivariant map $f$ induces a holomorphic map $f_S:S\fd\B^1$ such that
$f(w+2\pi b)=e^{i\theta}f(w)$ for every $w\in S$. But the holomorphic map
$g_S:S\fd\C$, $w\longmapsto e^{-i\frac{\theta}{2\pi b}w}f_S(w)$ descends to
$A_b$ and, since $|g_S|\leq e^{\frac{\theta}{2b}}|f_S|$, the previous arguments
applied to $g_S$ instead of $f_U$  do imply the same contradiction.
\end{demo}

The conclusion of
Theorem~\ref{surfaces} now easily follows. Indeed, denoting by $\om'_1$ the
K{\"a}hler form of the induced complete hyperbolic metric of finite volume
on $M\simeq \Gamma'\backslash\B^1$, we may apply 
the Schwarz-Pick lemma to obtain $0\leq\la f^\star\om_n,\om'_1\ra\leq 2$
pointwise. Moreover,
$$\tau(\rho)=\int_M f^\star\om_n=\frac{1}{2}\int_M
\la f^\star\om_n,\om'_1\ra\, dV'_1\leq \int_M dV'_1=-2\pi\chi(M)~.$$
Therefore $\la f^\star\om_n,\om'_1\ra = 2$ everywhere and
hence, if $\rho':\Gamma'\fd\PU(1,1)$ is the representation induced by $\rho$,
$f$ defines a $\rho'$-equivariant isometry from $\B^1$ onto $\B^1$. In
particular, $\rho'$ (and so $\rho$)
is injective, and $\rho'(\Gamma')=\rho(\Gamma)$ must be discrete. Finally, since
$f$ descends to an isometry between $\G'\backslash\B^1$ and $\rho(\G)\backslash
\B^1$, $\rho$ is a uniformization representation and $\rho(\G)$ is a
lattice in $\PU(1,1)$. 
\end{demo}

\subsection{The general case}\label{gc}\hfill

Let now $\Gamma$ be the fundamental group of a $p$-times punctured closed
orientable surface of negative Euler characteristic $M$ and $\rho:\Gamma\fd\PU(n,1)$
be a reductive representation. In the following, 
we shall prove Theorem~\ref{surfaces} in this general setting. 

As usual, we identify $\G$ with  
a non-uniform torsion-free lattice in $\PU(1,1)$ and $M$ with
$\G\backslash\B^1$, and we call $C_1,\ldots,C_p$ the cusps of $M$.  
We can assume that $\rho$ is not tame, namely that for some $1\leq k\leq
p$, $\rho$ maps the
peripheral elements corresponding to the punctures $m_1,\ldots,m_k$ to
hyperbolic isometries of $\PU(n,1)$.  

We begin by showing that the proof given in the 
previous section, until Claim~\ref{localdiffeo} included, still holds here.   
 
If one looks carefully at the preceding arguments, one sees that what is
really needed is the existence of a $\rho$-equivariant harmonic map
$f:\B^1\fd\Bn$ such that:

\begin{enumerate}
\item[(a)] the norm $\|\nab\df\|$ is square integrable and hence
  $f^\star\omega_n$ is an $L^2$ form;
\item[(b)] we have $\tau(\rho)=\frac{1}{2}\int_M\la
  f^\star\omega_n,\om_1\ra\,dV_1$, which together with (a) implies as
  before the inequality $|\tau(\rho)|\leq{\rm Vol}(M)$; 
\item[(c)] Claim~\ref{gaussbonnet} holds and hence the equality 
  $|\tau(\rho)|={\rm Vol}(M)$ implies as for tame representations that $f$ is 
  a $\rho$-equivariant immersion from $\B^1$ into a totally geodesic copy
  of $\B^1$ in $\Bn$.
\end{enumerate} 

According to Proposition~\ref{finenergy},
there is no finite energy $\rho$-equivariant
map $\B^1\fd\Bn$, so that we can not apply Theorem~\ref{existharmonic}. 
Nevertheless, we shall prove that:

\begin{enumerate}
\item[(i)]  there exists a (infinite energy) $\rho$-equivariant
  harmonic map $f:\B^1\fd\Bn$; 
\item[(ii)] we have a control on the energy density of $f$ at infinity;
\item[(iii)] this control implies (a), (b) and (c).
\end{enumerate}

We shall use the notations of the proofs of the tame
case.

\medskip

(i) Our method is the same as in \cite{JZ97}. We first fix a particular
$\rho$-equivariant map $\Phi:\B^1\fd\B^n$. As
usual, we only define it in each cusp $C_i$ of $M$ and extend it on the compact part
of $M$. We call $\g_i$ a generator of the fundamental group of
$C_i$ (when working on a single cusp, we shall drop the subscript $i$). The
notations are as in the proof of Proposition~\ref{finenergy}. Let $C$ be a cusp
of $M$. If $\rho(\g)$ is parabolic or elliptic, then we define  
$\Phi$ as we defined $f$ in Proposition~\ref{finenergy}. If $\rho(\g)$ is
hyperbolic, let 
$\varphi:[0,a]\fd\B^n$ be a map sending $[0,a]$ proportionally to arclength into
the axis of $\rho(\g)$ (which is a geodesic in $\B^n$), and such that
$\varphi(a)=\rho(\g)\varphi(0)$. Then, define $\Phi$ by
$\Phi(v,t)=\varphi(v)$ (we have $\|\d\Phi\|_{(v,t)}=e^{t}\frac{\delta}{a}$
where $\delta$ is the translation length of $\rho(\g)$).

For each $s\in\R_+$, let $M_s$ be the compact Riemann surface with boundary,
obtained from $M$ by deleting the end $\{t>s\}$ in each cusp of $M$, and let
$F_s$ be the restriction of the fiber bundle $F:=\B^1\times_\rho\B^n$ to
$M_s$. From \cite{Corlette92}, we get the existence of a harmonic section $f_s$ of
$F_s$ which agrees with $\Phi$ on $\partial M_s$ and whose energy satisfies
$E(f_s)\leq E(\Phi_{|M_s})$.

We want to prove that there
exists a strictly increasing sequence $(s_n)_{n\in\N}$ going to infinity, such
that $f_{s_n}$ converges uniformly on every compact subset of $M$. The limiting
function $f$ will be the (infinite energy) harmonic section of $F$ we are
looking for. Using the 
same argument as in \cite{Corlette92}, it is sufficient to prove that for each
$0\leq t\leq s$, the energy $E({f_s}_{|M_t})$ is bounded independently of
$s$.

Let $C$ be a cusp of $M$ such that $\rho(\gamma)$ is hyperbolic. If we set
$C_s=C\cap M_s$, from the proof of Proposition~\ref{finenergy}, it is immediate that

$$E({f_s}_{|C_s\backslash C_t})\geq E(\Phi_{|C_s\backslash C_t})=\frac{\delta^2}{2a} (e^s-e^t)~.$$

Moreover, denoting by $M'$ the union of the compact part of $M$ with the cusps
$C_{k+1},\dots,C_p$, and letting $M''={M\backslash M'}$
 ($M'_s=M'\cap M_s$, $M''_s=M''\cap M_s$), we have

$$E(\Phi_{|M'_s})\leq E(\Phi_{|M'})<+\infty$$
and so, because of the energy minimizing property of $f_s$,

$$E({f_s}_{|M_t})\leq E(f_s)-E(\Phi_{|M''_s\backslash M''_t})\leq
E(\Phi_{|M'\cup M''_t})$$
which is independent of $s$.

\medskip

(ii) For each $i\in\{1,\dots,k\}$,
we now define
$$\begin{array}{rcl}
\alpha_i:[0,+\infty) & \longrightarrow &\R \\
t & \longmapsto & \displaystyle\frac{1}{2}\int_0^{a_i} \bigl(\|\df\|^2_{(v,t)}-\|\d\Phi\|^2_{(v,t)}\bigr)e^{-t} dv~,
\end{array}
$$
where $(v,t)\in[0,a_i]\times [0,+\infty)$ are the usual horospherical
coordinates in the cusp $C_i$. The $\alpha_i$'s are non-negative functions
(see the proof of Proposition~\ref{finenergy}). 

We claim that the energy density $e(f)$ is
controlled by $e(\Phi)$, namely:
$$
\int_0^{+\infty}\alpha_i(t)dt<+\infty \mbox{ for any }i\in\{1,\dots,k\},\mbox{ and
  }E(f_{|M'})<+\infty~.
$$ 
Otherwise, there exists $r>0$ such
that
$$E\bigl(f_{|M_r}\bigr)\geq 2+E\bigl(\Phi_{|M'\cup M''_r}\bigr)<+\infty~.$$
We may choose $s>r$ such that
$$\bigl|e(f_s)-e(f)\bigr|\leq \frac{1}{{\rm Vol}(M_r)}$$
on $M_r$ and so,
$$E\bigl({f_s}_{|M_r}\bigr)\geq 1+E\bigl(\Phi_{|M'\cup M''_r}\bigr)~.$$
Now, we have
$$\begin{array}{rcl}
\displaystyle E(f_s)=\int_{M_s} e(f_s)\,dV_1 & \geq & 1+E\bigl(\Phi_{|M'\cup
  M''_r}\bigr)+E\bigl(\Phi_{|M''_s\backslash M''_r}\bigr)\\
& \geq &1+E\bigl(\Phi_{|M_s}\bigr)~.
\end{array}
$$
But $f_s$ and $\Phi$ coincide on $\partial M_s$ and this contradicts the energy
minimizing property of the harmonic map $f_s$.

\medskip

(iii) {\it Proof of (a).} First, we remark that the map $\Phi$ used to construct
$f$ is a rank one
harmonic map on $M''$, and in fact, it is totally geodesic. From the formula of
Eells-Sampson (see Lemma~\ref{nabladfl2}), we thus get $\Delta e(\Phi)=2e(\Phi)$
on $M''$. Let now be $e_\Phi:M\fd\R$ be a ($C^2$) function such that
${e_\Phi}_{|M'}=0$ and ${e_\Phi}_{|M''\backslash M''_1}=e(\Phi)_{|M''\backslash
  M''_1}$. If $K$ is the compact subset $\overline{M''_1\backslash M''_0}$, we have
$$ \Delta\bigl(e(f)-e_\Phi\bigr)=-\|\nabla\df\|^2+{\rm Scal}(f^\star
R^n)+2\bigl(e(f)-e_\Phi\bigr)$$
on $M\backslash K$.
Denoting by $\eta_R$ ($R>1$) the usual cut-off functions, and using Green's
formula, we obtain
$$\begin{array}{rcl}
\displaystyle\int_{M}\eta_R\Delta \bigl(e(f)-e_\Phi\bigr)\,dV_1
&=&\displaystyle \int_{M}\bigl(e(f)-e_\Phi\bigr)(\Delta\eta_R)\,
dV_1\\
&=& \displaystyle\int_{M'} e(f) (\Delta\eta_R)\,dV_1+\sum_{i=1}^k\int_0^{+\infty}(\Delta\eta_R)\alpha_i(t)dt+A~,
\end{array}
$$
because $\eta_R$ only depends on $t$ in the cusps (here $A$ is a constant
independent of $R$). Like in Lemma~\ref{nabladfl2}, the inequality
$\|\nabla\df\|^2\leq-\Delta\bigl(e(f)-e_\Phi\bigr) +2\bigl(e(f)-e_\Phi\bigr)$,
which is valid on $M\backslash K$,
now implies that $\|\nabla\df\|\in L^2(M)$ and, like in Lemma~\ref{canusef}, we
also obtain that $f^\star\om_n$ is $L^2$.\fin

\medskip

{\it Proof of (b).} Since $f^\star\om_n$ is $L^2$, we see from the proof of
Proposition~\ref{invharmonic} that we only need to show that each $f^\star
\varsigma_i$ is an $L^2$ form to get $\tau(\rho)=\frac{1}{2}\int_M \la
f^\star\omega_n,\omega_1\ra\,dV_1$. 

Let $i\in\{1,\dots,k\}$. We take
horospherical coordinates $(z',v',t')\in\B^n$ such that the axis of
$\rho(\gamma_i)$ is given by $z'=v'=0$. Then
$\Phi^\star\varsigma_i=\Phi^\star dz'=0$, that is $\|\d\Phi\|=\|\Phi^\star dt'\|$.
Since $\int_0^{a_i}\bigl(\|f^\star dt'\|_{(v,t)}^2-\|\Phi^\star dt'\|^2_{(v,t)}\bigr) dv\geq 0$
for each $t$, the fact that $\|\d f\|^2\geq\|f^\star d
t'\|^2+\|f^\star\varsigma_i\|^2$ together with (ii) implies that
$f^\star\varsigma_i$ is $L^2$.

The proof of the inequality $|\tau(\rho)|\leq{\rm Vol}(M)$ is like
in section~\ref{proof4}, the key point being that $\|\nabla\df\|\in L^2(M)$.\fin

\medskip

{\it Proof of (c).} We only have to show that for each $\varepsilon>0$,
$$
\lim_{r\pfd +\infty} \frac{{\rm Area}_{\phi^\star
    h_\varepsilon}\bigl(B_{\phi^\star
    h_\varepsilon}(r)\bigr)}{r^2}=0~.
$$ 
In a cusp $C_i$ for which $\rho(\g_i)$ is hyperbolic, we first compare
the volume element $dV_{\phi^\star h_\varepsilon}$ to $dV_{\e g_1+\Phi^\star
  g_n}=\sqrt{\e^2+2\e\, e(\Phi)} \,dV_1$ induced by the model map $\Phi$. For any
$s>0$, we
immediately get
$${\rm Area}_{\phi^\star h_\varepsilon}(C_{i,s})\leq {\rm Area}_{\e g_1+\Phi^\star
  g_n}(C_{i,s})+2\int_0^s\alpha_i(t)dt+\frac{1}{2\e}\int_{C_{i,s}}\|f^\star\om_n\|^2 dV_1~.$$
Recall that $f^\star\om_n$ is $L^2$ and, since $E(f_{|M'})<+\infty$, the
other cusps have finite area w.r.t. $\phi^\star h_\e$. So
$$\begin{array}{rcl}
\displaystyle\lim_{r\pfd +\infty} \frac{{\rm Area}_{\phi^\star
    h_\varepsilon}\bigl(B_{\phi^\star
    h_\varepsilon}(r)\bigr)}{r^2} &\leq &
\displaystyle\lim_{r\pfd +\infty} \frac{{\rm Area}_{\e g_1+\Phi^\star
  g_n}\bigl(B_{\phi^\star
    h_\varepsilon}(r)\cap M''\bigr)}{r^2}\\
&\leq& 
\displaystyle\lim_{r\pfd +\infty} \frac{{\rm Area}_{\e g_1+\Phi^\star
  g_n}\Bigl(B_{g_1}\bigl(\frac{r}{\sqrt\e}\bigr)\cap M''\Bigr)}{r^2}
\end{array}$$
where $B_{g_1}(r)$ is the geodesic ball of radius $r$ with respect to the metric
$g_1$.

But $B_{g_1}\bigl(\frac{r}{\sqrt\e}\bigr)\cap C_i\subset C_{i,{\frac{r}{\sqrt\e}}}$,
and $dV_{\e g_1+\Phi^\star
  g_n}=\sqrt{\e^2+\e\,\frac{\delta_i^2}{a^2_i}e^{2t}}\,e^{-t}dv\,dt$ on $C_i$. From this,
we conclude that there exists a positive constant $A_\varepsilon$ such that
$${\rm Area}_{\e g_1+\Phi^\star
  g_n}\Bigl(B_{g_1}\bigl(\frac{r}{\sqrt\e}\bigr)\cap M''\Bigr)\leq A_\e\, r~.$$
Thus the area of balls for the metric $\phi^\star h_\e$ grows at most
linearly and the result follows. \fin

\medskip

As in the tame case,
we therefore know that $\rho(\G)$ stabilizes a totally geodesic copy of
$\B^1$ in $\Bn$, and hence we can consider $\rho$ as a homomorphism from
$\G$ to $\PU(1,1)$. Moreover, the harmonic map $f:\B^1\fd\B^1$ is a local
diffeomorphism and therefore we can pull-back the complex structure of the
target to $M$. The uniformization theorem implies that there exists a
complete hyperbolic metric $g_1'$ on $M$ compatible with this new complex
structure and a discrete torsion-free subgroup $\G'=u(\G)$ of $\PU(1,1)$,
isomorphic to $\G$, such that $(M,g_1')$ is isometric to
$\G'\backslash\B^1$. 

However, contrary to the case of tame representations,
Claim~\ref{lattice} does not a priori hold and we will assume that 
$\G'$ is not a lattice, namely that the volume of $M$ with respect to
$g_1'$ is infinite (in fact, Claim~\ref{lattice} and the proof of
Claim~\ref{tame} below imply that $\G'$ is a lattice if and only if $\rho$ is
tame).  
This means that some punctures, say $m_1,\ldots,m_q$, have neighbourhoods of
infinite $g_1'$-volume. More precisely, peripheral elements of $\G$
corresponding to these punctures are sent by the new uniformization
representation $u$ to hyperbolic isometries of $\PU(1,1)$ and for each
$i\in\{1,\ldots,q\}$, $m_i$ has a neighbourhood $g_1'$-isometric to the
annulus 
$$
A_{i}:=\{z\in\C\, ,\,e^{-\pi/b_i}<|z|<a_i\},~b_i>0,~1>a_i>e^{-\pi/2b_i},
$$
endowed with the metric 
$$
{\Big(\frac{2b_i}{|z|\sin(b_i\log|z|)}\Big)}^2\,dzd\ov z~.
$$ 
In particular, there exist disjoint simple closed $g_1'$-geodesics 
$c_1,\ldots,c_q$ in $M$, corresponding to the circles
$\{|z|=e^{-\pi/2b_i}\}$ in $A_i$, such that $c_i$ is null-homotopic in 
$M\cup\{m_i\}$ for each $i\in\{1,\ldots,q\}$. 

If we cut $M$ along these geodesics and double the remaining finite volume
part $M_0$ (whose interior is diffeomorphic to $M$) along its (possibly
disconnected) $g_1'$-geodesic boundary, we obtain a surface $2M$,
on which $g_1'$ extends by symmetry to a complete hyperbolic metric of finite
volume. We call $\sigma$ the $g_1'$-isometric (antiholomorphic) involution of
$2M$. Note that the Euler characteristic $\chi(2M)$ of $2M$ equals
$2\chi(M)$.     

\medskip

We want to extend the representation $\rho:\G\fd\PU(1,1)$ to a
representation $2\rho$ of the fundamental group $2\G$ of $2M$. For this we
need an adapted presentation of $2\G$. 

The hyperbolic metric $g_1'$ on $2M$ allows us to identify $2\G$ with a
lattice in $\PU(1,1)$ and $2M$ with $2\G\backslash\B^1$. Call again
$\pi:\B^1\fd 2M$ the covering projection and  choose a connected component $X'$ of
$\B^1\backslash\bigcup_{i=1}^q\pi^{-1}(c_i)$. Since $\bigcup_{i=1}^q\pi^{-1}(c_i)$
consists of disjoint geodesics of $\B^1$, $X'$ is a convex set whose
stabilizer in $2\G$ can be identified with $\G'$. Moreover,
$\pi(X')=\G'\backslash X'$ is just the interior of $M_0$, seen as a subset of
$2M$. For each $i$ choose a lift $\tilde c_i$ of $c_i$ in the boundary of
$X'$, a generator $\g_i'$ of the cyclic stabilizer of $\tilde c_i$ in
$\G'$, and call $\sigma_i$ the symmetry w.r.t. $\tilde c_i$. It is
easy to see that the symmetries $\sigma_i$ are lifts of $\sigma$ and
therefore for all $i$ there exists $h_i\in 2\G$ such that $\sigma_i =
h_i\sigma_1$.   

Let $X''=\sigma_1 X'$ be the other component of
$\B^1\backslash\bigcup_{i=1}^q\pi^{-1}(c_i)$ adjacent to $\tilde c_1$, and call
$\G''$ ($=\sigma_1\G'\sigma_1$) its stabilizer in $2\G$.    
Note that, for all $i$, $h_i^{-1}\g_i' h_i
=\sigma_1\sigma_i\g_i'\sigma_i\sigma_1=\sigma_1\g_i'\sigma_1$ belongs to
$\G''$.

The fundamental group $2\G$ of $2M$ has the following abstract
presentation (see for example~\cite{Serre}). 
It is generated by the sets $\G'$ and $\G''$, together
with the $h_i$'s, subject to the relations:\\
- the relations of the groups $\G'$ and $\G''$;\\
- for all $i\in\{1,\ldots,q\}$, the element $\g_i'$ of $\G'$ is identified
with the element $\g_i'':=h_i^{-1}\g_i' h_i$ of $\G''$.  

\medskip

Now, it follows from the discussion in the proof of Claim~\ref{lattice}
that for all $i\in\{1,\ldots,q\}$, $g_i:=\rho\circ u^{-1}(\g_i')$ is an
hyperbolic isometry of $\B^1$. Call $s_i$ the symmetry w.r.t. the axis of
$g_i$ and set
$$
\left\{
\begin{array}{l}
2\rho(\g')=\rho\circ u^{-1}(\g')\mbox{ for } \g'\in\G'~; \\  
2\rho(\g'')=s_1(\rho\circ u^{-1}(\sigma_1\g''\sigma_1))s_1
\mbox{ for }\g''\in\G''~;\\
2\rho(h_i)=s_i s_1\mbox{ for } i\in\{1,\ldots,q\}~.
\end{array}
\right.
$$
This is clearly compatible with the relations of $2\G$ and hence $2\rho$ is
a well-defined representation of $2\G$ into $\PU(1,1)$. 

\begin{claim}\label{tame}
  The representation $2\rho:2\G\fd\PU(1,1)$ is tame.  
\end{claim}

\begin{demo}
It is enough to show that, for any $i\in\{q+1,\ldots,p\}$, the peripheral
elements  $\g_i$ corresponding to the punctures $m_i$ are not mapped by
$\rho$ to hyperbolic isometries of $\B^1$. 

Suppose that one of them is. We drop the subscript $i$. The quotient
$\la u(\g)\ra\backslash\B^1$ is holomorphically equivalent to a punctured disc
$D^\star=\{z\in{\mathbb C}\, ,\,0<|z|<1\}$, whereas the quotient
$\la\rho(\g)\ra\backslash\B^1$ is holomorphically equivalent to an annulus
$A=\{z\in{\mathbb C}\, ,\,a<|z|<1\}$. The $\rho$-equivariant holomorphic map
$f:\B^1\fd\B^1$ induces a holomorphic map $f$ from $D^\star$ to $A$. This
map is bounded and hence extends to a holomorphic map from the whole disc
$D$ to $\ov A$. This is impossible since the $\rho$-equivariance of $f$ implies
that for all $0<r<1$, the loops $f(\{z\in D^\star,\,|z|=r\})$ are homotopically
non trivial.  
\end{demo}

\begin{claim}
 The invariant $\tau(2\rho)$ is maximal: $\tau(2\rho)=-2\pi\chi(2M)$. 
\end{claim}

\begin{demo}
  Since $2\rho$ is tame, there exists a finite energy $2\rho$-equivariant map
  $f:\B^1\fd\B^1$ and $\tau(2\rho)=\int_{2M} f^\star\om_1$.   

For each $i\in\{1,\ldots,q\}$, we may choose a symmetric tubular neighbourhood  
$N_i$ of the geodesic $c_i$ and a lift $\wt N_i$ in $\B^1$ such that
$N_i=\la\g_i'\ra\backslash\wt N_i$. Let $\psi_i$ be a Kähler potential
associated as before to a fixed point $\xi_i$ of $2\rho(\g_i')$ and let
$\varsigma_i$ be the corresponding 1-form. Then 
$f^\star\varsigma_i$ is a 1-form on $N_i$ such that $f^\star\om_1=\d
f^\star\varsigma_i$. If now $\eta$ is a symmetric cut-off function on $2M$
identically equal to 0 outside the $N_i$'s and to 1 close to the geodesics
$c_i$, we get: 
$$
\begin{array}{rcl}
\tau(2\rho) &=& \ds\int_{2M} \Big[f^\star\om_1
-\d\Big(\sum_{i=1}^q\eta f^\star\varsigma_i\Big)\Big]\\
& = &
 \ds\int_{M_0} \Big[f^\star\om_1
-\d\Big(\sum_{i=1}^q\eta f^\star\varsigma_i\Big)\Big]
+ \ds\int_{\sigma M_0} \Big[f^\star\om_1
-\d\Big(\sum_{i=1}^q\eta f^\star\varsigma_i\Big)\Big]\\
& = &
 \ds\int_{M_0} \Big[f^\star\om_1
-\d\Big(\sum_{i=1}^q\eta f^\star\varsigma_i\Big)\Big]
- \ds\int_{M_0} \Big[\sigma^\star f^\star\om_1
-\d\Big(\sum_{i=1}^q\eta \sigma^\star f^\star\varsigma_i\Big)\Big]~.\\
\end{array}
$$

Extending the form $f^\star\om_1-\d(\sum_{i=1}^q\eta f^\star\varsigma_i)$ on
$M_0$ by zero to a form on $M$, one sees that $\int_M
[f^\star\om_1-\d(\sum_{i=1}^q\eta f^\star\varsigma_i)]=\tau(\rho)$.

The forms $\sigma^\star f^\star\om_1$ on $M$ and $\sigma^\star
f^\star\varsigma_i$ on $N_i$ are respectively induced by $\sigma_1^\star
f^\star\om_1$ on $\B^1$ and $\sigma_1^\star f^\star\varsigma_i$ on
$\sigma_1 \wt N_i$. Since $s_1^\star\om_1=-\om_1$, one has $\sigma_1^\star
f^\star\om_1=-(s_1\circ f\circ\sigma_1)^\star\om_1$. Moreover,  
$\varsigma_i=s_1^\star s_1^\star(-\d^c\psi_i)= s_1^\star\d^c(\psi_i\circ s_1)$, 
and $\psi_i\circ s_1=\psi_i\circ s_i\circ s_1=2\rho(h_i)^\star\psi_i$. Hence
$\sigma_1^\star f^\star\varsigma_i=-(s_1\circ f\circ\sigma_1)^\star
2\rho(h_i)^\star\varsigma_i$. Since $s_1\circ f\circ\sigma_1$ is also
$\rho$-equivariant, the 1-form $(s_1\circ f\circ\sigma_1)^\star
2\rho(h_i)^\star\varsigma_i=h_i^\star (s_1\circ f\circ\sigma_1)^\star\varsigma_i$ on
$h_i^{-1}\wt N_i=\sigma_1\wt N_i$ induces $(s_1\circ
f\circ\sigma_1)^\star\varsigma_i$ on $N_i$. Finally,
$$
 - \int_{M_0} \Big[\sigma^\star f^\star\om_1
-\d\Big(\sum_{i=1}^q\eta \sigma^\star f^\star\varsigma_i\Big)\Big]
 = 
 \int_{M_0}  \Big[(s_1\circ f\circ\sigma_1)^\star\om_1
-\d\Big(\sum_{i=1}^q \eta (s_1\circ f\circ\sigma_1)^\star\varsigma_i\Big)\Big]\\
$$      
and the r.h.s. again equals $\tau(\rho)$ since $s_1\circ f\circ\sigma_1$ is a
$\rho$-equivariant map.  
Hence $\tau(2\rho)=2\tau(\rho)=-4\pi\chi(M)=-2\pi\chi(2M)$ and the lemma is proved.
\end{demo}

The representation $2\rho:2\G\fd\PU(1,1)$ is therefore a tame representation of
maximal invariant. It follows from the results of section~\ref{wtpr} that
there exists a $2\rho$-equivariant
isometry $f$ from $\B^1$ onto $\B^1$. Since $f$ is a fortiori $(\rho\circ
u^{-1})$-equivariant, $\rho$ is a uniformization
representation and we are done.

\end{document}